\documentclass[11pt,a4paper]{article}

 \usepackage{geometry}
\newgeometry{ hmargin={20mm,20mm}}

\usepackage{algorithm}
\usepackage{amsthm}
\usepackage{amsmath}
\usepackage{amssymb}
\usepackage{bm}
\usepackage{booktabs}
\usepackage{caption,subcaption}

\usepackage[numbers]{natbib}

\usepackage{tikz}
\usetikzlibrary{arrows.meta,
                chains,
                positioning,
                shapes.geometric
                }

\DeclareMathOperator*{\argmin}{argmin}
\DeclareMathOperator*{\minimize}{minimize}
\DeclareMathOperator*{\maximize}{maximize}

\DeclareMathOperator{\diag}{diag}

\newtheorem{theorem}{Theorem}
\newtheorem{cor}{Corollary}

\newtheorem{prop}{Proposition}
\newtheorem{defn}{Definition}

\DeclareMathOperator{\bone}{\bf 1}
\DeclareMathOperator{\bzero}{\bf 0}
\DeclareMathOperator{\bA}{\bf A}
\DeclareMathOperator{\ba}{\bf a}
\DeclareMathOperator{\bB}{\bf B}
\DeclareMathOperator{\blb}{\bf b}
\DeclareMathOperator{\bc}{\bf c}
\DeclareMathOperator{\bdA}{\bf dA}
\DeclareMathOperator{\bdb}{\bf db}
\DeclareMathOperator{\bdG}{\bf dG}
\DeclareMathOperator{\bdh}{\bf dh}
\DeclareMathOperator{\bdp}{\bf dp}
\DeclareMathOperator{\bdQ}{\bf dQ}
\DeclareMathOperator{\bdz}{\bf dz}

\DeclareMathOperator{\bd}{\bf d}

\DeclareMathOperator{\bF}{\bf F}
\DeclareMathOperator{\bG}{\bf G}

\DeclareMathOperator{\bI}{\bf I}
\DeclareMathOperator{\bh}{\bf h}
\DeclareMathOperator{\bl}{\bf l}
\DeclareMathOperator{\bM}{\bf M}
\DeclareMathOperator{\bp}{\bf p}
\DeclareMathOperator{\bQ}{\bf Q}

\DeclareMathOperator{\br}{\bf r}

\DeclareMathOperator{\bs}{\bf s}
\DeclareMathOperator{\bU}{\bf U}
\DeclareMathOperator{\bu}{\bf u}

\DeclareMathOperator{\bv}{\bf v}
\DeclareMathOperator{\bW}{\bf W}
\DeclareMathOperator{\bw}{\bf w}

\DeclareMathOperator{\bx}{\bf x}

\DeclareMathOperator{\by}{\bf y}
\DeclareMathOperator{\bz}{\bf z}

\DeclareMathOperator{\bdlambda}{ \bd \boldsymbol{\lambda} }
\DeclareMathOperator{\bdnu}{ \bd \boldsymbol{\nu}}
\DeclareMathOperator{\bepsilon}{\bm \epsilon}

\DeclareMathOperator{\betta}{\bm \eta}

\DeclareMathOperator{\blambda}{\bm \lambda}
\DeclareMathOperator{\bmu}{\bm \mu}
\DeclareMathOperator{\bnu}{\bm \nu}

\DeclareMathOperator{\btheta}{\bm \theta}

\providecommand{\keywords}[1]{\textbf{\textit{Keywords:}} #1}

\newcommand{\toi}[2][i]{%
  \mathop{
    \mathrm{#2}^{( #1 )}
  }
}

\newcommand{\toit}[2][i]{%
  \mathop{
    \mathrm{#2}^{{T}^{( #1 )}}
  }
}

\begin{document}


\title{Efficient differentiable quadratic programming layers: an ADMM approach}
\author{Andrew Butler and Roy H. Kwon \\ University of Toronto\\Department of Mechanical and Industrial Engineering}
\maketitle

\begin{abstract}

Recent advances in neural-network architecture allow for seamless integration of convex optimization problems as differentiable layers in an end-to-end trainable neural network. Integrating medium and large scale quadratic programs into a deep neural network architecture, however, is challenging as solving quadratic programs exactly by interior-point methods has worst-case cubic complexity in the number of variables. In this paper, we present an alternative network layer architecture based on the alternating direction method of multipliers (ADMM) that is capable of scaling to problems with a moderately large number of variables. Backward differentiation is performed by implicit differentiation of the residual map of a modified fixed-point iteration. Simulated results demonstrate the computational advantage of the ADMM layer, which for medium scaled problems is approximately an order of magnitude faster than the OptNet quadratic programming layer. Furthermore, our novel backward-pass routine is efficient, from both a memory and computation standpoint, in comparison to the standard approach based on unrolled differentiation or implicit differentiation of the  KKT optimality conditions.  We conclude with examples from portfolio optimization in the integrated prediction and optimization paradigm.

\end{abstract}


\section{Introduction}\label{sec:intro}

Many problems in engineering, statistics and machine learning require solving convex optimization programs. In many real-world applications, the solution to the convex optimization program is a single component in a larger decision-making process (see for example \citep{Butler2020IPO,Butler2021IPOb,Donti2017}). Recent advances in neural-network architecture allow for seamless integration of convex optimization programs as differentiable layers in an end-to-end trainable neural network \citep{Agrawal2019,Agrawal2020,Amos2017}.

In this paper, we consider convex programming layers that take the form of parametric quadratic programs (PQPs) with linear equality and box inequality constraints:
\begin{equation}\label{eq:pqp}
\begin{split}
\minimize_{\bz} \quad &   \frac{1}{2} \bz^T \bQ(\btheta)  \bz + \bz^T \bp(\btheta)\\
\text{subject to} \quad &   \bA(\btheta)  \bz = \blb(\btheta), \quad \bl(\btheta) \leq \bz \leq \bu(\btheta) \\
\end{split}
\end{equation}
Here $\bz \in \mathbb{R}^{d_z}$ denotes the decision variable and $\bQ(\btheta)$, $\bp(\btheta)$, $\bA(\btheta)$, $\blb(\btheta)$, $\bl(\btheta)$, $\bu(\btheta)$ are the parameterized problem variables.   Program $\eqref{eq:pqp}$ occurs in many applications to statistics \citep{Tibshirani1996, Tikhonov1963}, machine-learning \citep{Ganti2011,Ho2015}, signal-processing \citep{Kim2008} and finance \citep{ Black1991, Michaud2008b,Markowitz1952}.  

In general, a differentiable convex optimization layer can be viewed as a function that maps the program input variables to optimal primal(-dual) solution(s). Therefore, in a fully integrated system, optimizing problem variables by  backpropogation ultimately requires computing the action of the Jacobian of the optimal solution(s) with respect to the corresponding program input variables. For example, the OptNet layer, proposed by \citet{Amos2017} is a specialized differentiable optimization layer that uses a primal-dual interior-point method for small scale batch quadratic programs. The authors demonstrate that the solution to the system of equations provided by the KKT conditions at optimality provide a system of equations for implicit differentiation with respect to all relevant problem variables. Therefore,  by strategically caching the KKT matrix factorization(s) (obtained during the forward-pass), then the gradients required for backpropogation are obtained at little extra additional cost. Indeed, the authors note that the OptNet layer is computationally efficient and therefore practical for small-scale problems $(d_z < 100)$. However, solving convex quadratic programs exactly by interior-point methods has worst-case time complexity on the order of $\mathcal{O}(d_z^3)$  \citep{Goldfarb1991}. Therefore, for medium $(100 < d_z < 1000)$ and large  $(d_z > 1000)$ scale problems, embedding an OptNet layer in a larger neural network can be computationally intractable.

In this paper, we address the computational challenges in the medium to large scale limit and propose an alternative differentiable network architecture for batch constrained quadratic programs of the form of Program $\eqref{eq:pqp}$. Our differentiable quadratic programming layer is built on top of the alternating direction method of multipliers (ADMM) algorithm, which recently has become increasingly popular for solving large-scale convex optimization problems   \citep{Boyd2011, Boyd2016, Sopa2019, Stellato2020, Schu2020}. Indeed, embedding the ADMM algorithm in a larger neural-network system has been fundamental to recent innovations in signal processing, compressed sensing, imaging and statistics (see for example \citep{ Diamond2018, Xie2019, Yang2017}). The standard ADMM-network implementation computes the relevant gradients by unrolling the ADMM computational graph, which is memory inefficient and typically requires substantially larger networks. Furthermore, when the number of ADMM iterations is large in the forward-pass then the unrolled gradient can be computationally demanding. Alternatively, many differentiable convex optimization layers compute the relevant gradients by implicitly differentiating the Karush–Kuhn–Tucker (KKT) optimality conditions. Indeed differentiating the KKT system of equations is possible but unfortunately requires solving a system of equations of dimension $\mathbb{R}^{3d_z \times 3d_z}$, which can also be computationally impractical. As an alternative, we present a novel modified backward-pass routine that is efficient from both a memory and computational standpoint. Specifically, we recast the ADMM algorithm as a fixed-point iteration and apply the implicit function theorem to the resulting residual map in order to obtain the relevant backpropogation gradients. A primary advantage of our fixed-point differentiation is that the fixed-point iteration is of dimension $d_z$ and therefore the resulting system of equations is approximately $3$ times smaller than the equivalent KKT system.  Finally, our differentiable ADMM-layer and all algorithmic implementations is made available as an open-source R package, available here:
$$ \text{https://github.com/adsb85/lqp}$$

The remainder of the paper is outlined as follows. We begin with the problem motivation and a brief discussion of related work in the field of differentiable convex optimization layers. In Section \ref{sec:method_admm} we review the ADMM algorithm and present our ADMM-based architecture for forward solving batch PQPs. We review the KKT based implicit differentiation and then present the fixed-point iteration and derive the  expression for the relevant gradients. In  Section \ref{sec:results} we perform several simulations (under numerous model specifications) and compare the computational efficiency and performance accuracy of our ADMM-layer with the state-of-the-art OptNet layer.  We demonstrate that for medium-scale problems, our  ADMM-layer implementation is  approximately an order of magnitude faster than the OptNet layer and provides solutions that are equally as optimal. Moreover, we compare the computational efficiency of the backpropogation routines based on unrolled differentiation, KKT implicit differentiation and fixed-point implicit differentiation. We demonstrate  that the fixed-point implicit differentiation is universally more efficient than the KKT implicit differentiation, and under certain conditions is preferred to the unrolled differentiation.  We conclude with a real world application of the ADMM-layer to a medium scale portfolio optimization problem in the integrated prediction and optimization paradigm.

\subsection{Related work:}

\subsubsection{Problem Motivation}

Many problems in machine learning, statistics, engineering and operations research involve both predictive forecasting and decision based optimization. Recent advances in neural network architecture embed convex optimization programs as differentiable layers in a larger neural network structure. A fully integrated prediction and optimization (IPO) architecture therefore enables the integration of predictive forecasting and convex optimization and allows for the minimization of the total decision error induced by the forecast estimates (see for example \citep{Butler2021IPOb,Butler2020IPO,Donti2017,Elma2020,Elma2020b}). This is in contrast to a more traditional `predict, then optimize' approach which  would first fit the predictive models (for example by maximum likelihood or least-squares) and then `plug-in' those predictions to the corresponding decision-based optimization program. While it is true that a perfect predictive model would lead to optimal decision making, in reality, all predictive models do make some error, and thus an inefficiency exists in the `predict, then optimize' paradigm. With the widespread adoption of machine-learning and data science in operations research, there has been a growing body of literature on data-driven decision making and the relative merits of decoupled versus integrated predictive decision-making \citep{Amos2019, Bert2020, Grigas2021, Mandi2019, Mandi2020, Uysal2021}.

The preliminary findings of the aforementioned work advocate strongly for an IPO approach. Indeed, IPO models typically exhibit lower model complexity and demonstrate statistically significant improvements in out-of-sample performance in comparison to more traditional `predict, then optimize' models. Unfortunately, training medium and large scale IPO models can be computationally demanding. For example, \citet{Elma2020, Elma2020b} report that their `smart predict, then optimize' models can take several hours to train medium and large scale problems, while traditional prediction methods typically take seconds or minutes to train. \citet{Butler2021IPOb} provide details on the computational complexity of  IPO models trained with the OptNet layer and demonstrate that when $d_z > 100$ then fitting IPO models can be computationally impractical. Improving the efficiency of the integrated framework is therefore an open and important area of research.

\subsubsection{Differentiable convex optimization layers:}

Differentiable convex optimization layers provide an efficient and seamless framework for integrating predictive forecasting with downstream decision-based optimization in a end-to-end trainable neural network. Modern neural network technology (such as \textit{torch} or \textit{tensorflow}) require that every layer in the network inherits a forward and backward-pass routine. For differentiable convex optimization layers , the forward-pass is typically an iterative optimization algorithm that converts problem variables into optimal primal(-dual) solution(s). The backward-pass routine therefore computes the action of the Jacobian of the optimal solution(s) with respect to all problem variables, and returns the left matrix-vector product of the Jacobian with the previous backward-pass gradient(s).

For example, the state-of-the-art OptNet layer implements a primal-dual interior-point method for solving small-scale batch constrained quadratic programs \citep{Amos2017}. For backward differentiation, the authors implement a novel and efficient argmin differentiation routine that implicitly differentiates the KKT system of equations at optimality. By strategically caching the factorized KKT left-hand side matrix then the resulting method is shown to be computationally tractable for small problems within the context of deep neural network architectures. The author's acknowledge, however, that the OptNet layer may be impractical for  optimization problems with a moderate to large number of variables.

More recently,  \citet{Agrawal2020} provide a general framework for differentiable convex cone programming. Their forward-pass recasts the conic program in its equivalent homogeneous self-dual embedding form, which is then solved by operator splitting \citep{Boyd2016}.  In the backward-pass, the relevant gradients are obtained by implicit differentiation of the residual map provided by the homogeneous self-dual embedding. The resulting differentiable cone programming layer is flexible, but requires the user to transform their problem into a canonical form, which is often time-consuming, prone to error and requires a certain level of domain expertise.

Alternatively, \citet{Agrawal2019}, provide a domain-specific language for differentiable disciplined convex programs. Their approach abstracts away the process of converting problems to canonical form with minimal loss in computational efficiency in comparison to specialized convex optimization layers. They also provide an efficient sparse matrix solver, which for sparse quadratic programs is on average an order of magnitude faster than the OptNet layer. Similarly, \citet{Blondel2021} provide an efficient and modular approach for implicit differentiation of optimization problems. They consider KKT, proximal gradient and mirror descent fixed-point implicit differentiation  and provide a software infrastructure for efficiently integrating their modular implicit differentiation routines with state-of-the-art optimization solvers. That said, for solving batches of convex optimization problems it is often preferred and more efficient to avail of optimization solvers that have the ability to exploit fast GPU-based batch solves.

\subsubsection{ADMM and unrolled differentiation:}

The alternating direction method of multipliers (ADMM) algorithm, first proposed by \citet{Gabay1976} and \citet{Glow1975}, is well suited to many large-scale and distributed problems common to applications of statistics, machine learning, control and finance. We note that the ADMM algorithm is closely related to algorithms such as dual ascent, the augmented Lagrangian method of multipliers, and operator (Douglas–Rachford) splitting and refer to \citet{Boyd2011} for a comprehensive overview.

Embedding the above mentioned algorithms in larger neural-network structures has been fundamental to recent innovations in signal processing, compressed sensing, imaging and statistics (see for example \citep{Belanger2017, Diamond2018, Domke2012, Feng2017, Lorraine2018}). Perhaps most closely related to our own work, the ADMM-Net, first proposed by \citet{Yang2017}, recasts and embeds the iterative ADMM procedure as a fully learnable network graph. The authors provide examples from compressive sensing magnetic resonance imaging and demonstrate that their ADMM-Net achieves state-of-the-art model accuracy and computationally efficiency. More recently, the Differentiable Linearized ADMM (D-LADMM) algorithm, proposed by  \citep{Xie2019}, generalizes the ADMM-Net and is capable of solving general deep learning problems with equality constraints. The authors show that there exists a set of learnable parameters for D-LADMM to generate global solutions and provide the relevant convergence analysis. However, in all cases mentioned above the authors consider either unconstrained or linear equality constrained least-squares problems and do not consider inequality constraints. Furthermore,  they perform the action of argmin differentiation by unrolling the `inner-loop' of the convex optimization routine, which  necessitates substantially larger and less efficient networks \citep{Amos2017, Amos2019}. 

Our ADMM-layer derives inspiration from both the differentiable convex optimization and ADMM network literature. To our knowledge, our ADMM-layer is the first implementation of its kind that can efficiently handle medium to large scale differentiable constrained PQPs. In this paper we demonstrate that solving medium and large scale QPs by interior-point methods can be  computationally burdensome. An obvious course of action is to replace the interior-point algorithm in the forward-pass with a more computationally efficient first-order method, such as ADMM. However, implementing an efficient backward-pass routine by implicit differentiation of  the KKT conditions is challenging as the ADMM algorithm does not explicit solve the KKT system of equations. Therefore, unlike the OptNet implementation, at each `outer' iteration (hereafter referred to as epochs) we must form and solve the resulting KKT system, which for large scale problems creates a computational bottleneck.  Unrolled differentiation of the ADMM algorithm is appropriate for small scale problems. However, for larger scale problems or for problems that require solving the convex optimization problem to a high degree of accuracy, an unrolled differentiation approach can also be computationally impractical. In contrast, our novel fixed-point implicit differentiation method is shown to be computationally efficient and invariant to the number of `inner' iterations performed in the ADMM forward-pass. Furthermore, as mentioned earlier, a primary advantage of the fixed-point differentiation is that the fixed-point iterative scheme is of dimension $d_z$ and is therefore approximately $3$ times smaller than the KKT system.  We demonstrate that in the absence of a pre-factorized KKT system, the fixed-point implicit differentiation is preferred to the KKT implicit differentiation and is competitive with the efficient KKT factorization caching provided in the OptNet layer. Furthermore, we demonstrate that for small-scale problems, our ADMM-layer is competitive with the state-of-the-art OptNet layer in terms of accuracy and computational efficiency. For medium and large scale problems, our ADMM layer is shown to be approximately an order of magnitude faster than the OptNet layer. Of course, our ADMM-layer is not without its own limitations and areas for improvement, discussed in detail in Section \ref{sec:conclusion}.

\section{Methodology}\label{sec:method_admm}
In general, the ADMM algorithm  is applied to problems of the form:
\begin{equation}\label{eq:admm}
\begin{split}
\minimize \quad &   f(\bx) + g(\bz) \\
\text{subject to} \quad &   \bA\bx + \bB\bz = \bc \\
\end{split}
\end{equation}
with decision variables $\bx \in \mathbb{R}^{d_x}$, $\bz \in \mathbb{R}^{d_z}$ and problem variables $\bA \in \mathbb{R}^{d_{eq} \times d_x}$, $\bB \in \mathbb{R}^{d_{eq} \times d_z}$ and $\bc \in \mathbb{R}^{d_{eq}}$. In order to guarantee convergence we assume that $f \colon \mathbb{R}^{d_x}  \to \mathbb{R}$ and  $g \colon \mathbb{R}^{d_z}  \to \mathbb{R}$ are closed, proper convex functions \citep{Boyd2011}. The augmented Lagrangian of Program $\eqref{eq:admm}$ is given by:
\begin{equation}\label{eq:admm_L}
L_{\rho}(\bx, \bz, \by) = f(\bx) + g(\bz) + \blambda^T( \br ) + \frac{\rho}{2}\lVert \br \rVert_2^2,
\end{equation}
with Lagrange dual variable $\blambda$, residual $\br = \bA\bx + \bB\bz - \bc$ and user-defined penalty parameter $\rho > 0$. We denote $\bmu = \rho^{-1}\blambda$ and therefore the well-known scaled ADMM iterations are as follows:

\begin{equation}\label{eq:admm_iter}
\begin{split}
\bx^{k+1} & = \argmin_{\bx} f(\bx) + \frac{\rho}{2} \lVert \bA\bx + \bB\bz^k - \bc + {\bmu}^k \rVert_2^2\\
\bz^{k+1} & = \argmin_{\bz} g(\bz) + \frac{\rho}{2} \lVert \bA\bx^{k+1} + \bB\bz - \bc + {\bmu}^k \rVert_2^2\\
{\bmu}^{k+1} & = {\bmu}^{k} +  \bA\bx^{k+1} + \bB\bz^{k+1} - \bc
\end{split}
\end{equation}
where $\bx^{k}$ and $\bz^{k}$ denote the decision variables at iteration $k$.

We denote $\br^k = \bA\bx^k + \bB\bz^k - \bc$ and $\bs^k = \rho \bA^T\bB(\bz^{k} - \bz^{k-1})$ as the primal and dual residual at iteration $k$. Let $\epsilon^p > 0$ and  $\epsilon^d > 0$ be the user defined stopping tolerances for the primal and dual residuals, respectively. Therefore a  reasonable stopping criteria would be:
\begin{equation}\label{eq:admm_stop}
\br^k \leq  \epsilon^p \quad \text{and} \quad \bs^k \leq  \epsilon^d.
\end{equation}

\subsection{ADMM for parametric quadratic programs}\label{sec:method_pqp_admm}
We consider convex parametric quadratic programs (PQPs) of the form:
\begin{equation}\label{eq:pqp_main}
\begin{split}
\minimize_{\bz} \quad &   \frac{1}{2} \bz^T \bQ(\btheta)  \bz + \bz^T \bp(\btheta) \\
\text{subject to} \quad &   \bA(\btheta)  \bz = \blb(\btheta), \quad  \bl(\btheta) \leq \bz \leq \bu(\btheta),
\end{split}
\end{equation}
with decision variable $\bz \in \mathbb{R}^{d_z}$. The objective function is therefore defined by a vector  $\bp(\btheta) \in  \mathbb{R}^{ d_z}$  and symmetric positive definite matrix $\bQ(\btheta) \in  \mathbb{R}^{ d_z \times d_x}$. Here, $\bA(\btheta) \in \mathbb{R}^{d_{eq} \times d_z}$, $\blb(\btheta)  \in  \mathbb{R}^{ d_{eq} }$ , $\bl(\btheta)  \in  \mathbb{R}^{ d_z}$ and $\bu(\btheta)  \in  \mathbb{R}^{ d_z}$ define the linear equality and box inequality constraints. We assume that all problem variables are parameterized by  $\btheta$ and are therefore trainable when integrated in an end-to-end neural network; rather than simply being supplied by the user.

\subsubsection{ADMM-layer: forward-pass}
Our ADMM-layer solves Program $\eqref{eq:pqp_main}$ in the forward-pass by applying the ADMM algorithm as outlined in Section \ref{sec:method_admm}. Note that for ease of notation we temporarily drop the parameterization, $\btheta$. 

We define
\begin{equation}\label{eq:f_admm}
f(\bx) = \frac{1}{2} \bx^T \bQ \bx + \bx^T \bp,
\end{equation}
with domain $\{\bx | \bA \bx = \blb\}$. Similarly, we define
\begin{equation}\label{eq:g_admm}
g(\bz) = \mathbb{I}_{\bl \leq \bz \leq \bu}(\bz) 
\end{equation}
where $\mathbb{I}_{\bl \leq \bz \leq \bu}(\bz)$ denotes the indicator function with respect to the linear inequality constraints. Program $\eqref{eq:pqp_main}$ is then recast to the following convex optimization Program:
\begin{equation}\label{eq:ppqp_admm}
\begin{split}
\minimize \quad &   f(\bx) + g(\bz)  \\
\text{subject to} \quad &   \bx - \bz = 0.
\end{split}
\end{equation}
Applying the ADMM iterations, as defined by Equations $\eqref{eq:admm_iter}$, to Program $\eqref{eq:ppqp_admm}$  therefore gives the following iterative optimization algorithm:
\begin{subequations} \label{eq:admm_qp_iter}
\begin{align}
\bx^{k+1} & = \argmin_{\{\bx |  \bA \bx = \blb\}} \frac{1}{2} \bx^T \bQ \bx + \bx^T \bp + \frac{\rho}{2} \lVert \bx - \bz^k + {\bmu}^k \rVert_2^2 \label{eq:admm_qp_x}\\
\bz^{k+1} & = \argmin_{\{ \bl \leq \bz \leq \bu \} } \frac{\rho}{2} \lVert \bx^{k+1} - \bz + {\bmu}^k \rVert_2^2  \label{eq:admm_qp_z}\\
{\bmu}^{k+1} & = {\bmu}^{k} +  \bx^{k+1} - \bz^{k+1} \label{eq:admm_qp_mu}
\end{align}
\end{subequations}

The ADMM algorithm, as defined by Equations $\eqref{eq:admm_qp_iter}$, allows for efficient optimization of medium and large scale quadratic programs. Firstly, we note that $\eqref{eq:admm_qp_z}$ is a least squares problem with box-inequality constraints, and therefore can be solved analytically. Specifically, we define the euclidean projection onto a set of box constraints as:
\begin{equation}\label{eq:project_box}
\Pi(\bx) = \begin{cases}
                \bl_j & \text{if } \bx_j < \bl_j\\
                \bx_j & \text{if } \bl_j \leq \bx_j \leq \bl_j\\
                \bu_j & \text{if } \bx_j > \bu_j\\
                \end{cases}.
\end{equation}
The analytic solution to Program $\eqref{eq:admm_qp_z}$ is therefore given by:
\begin{equation}\label{eq:admm_qp_z_simp}
\bz^{k+1}  = \Pi( \bx^{k+1}  + {\bmu}^k ).
\end{equation}
Furthermore, Program $\eqref{eq:admm_qp_x}$ is an equality constrained quadratic program, which can also be solved analytically. Specifically, the KKT optimality conditions of Program $\eqref{eq:admm_qp_x}$ can be expressed as the solution to the following linear system of equations:
\begin{equation} \label{eq:admm_qp_x_simp}
\begin{split}
\begin{bmatrix}
\bQ + \rho \bI_{\bx} &   \bA^T \\
\bA & 0
\end{bmatrix}
\begin{bmatrix}
\bx^{k+1}\\
\betta^{k+1}
\end{bmatrix}
= -
\begin{bmatrix}
 \bp - \rho (\bz^{k} - \bmu^{k})\\
 -\blb
\end{bmatrix},
\end{split}
\end{equation}
with identity matrix $\bI_{\bx} \in \mathbb{R}^{d_x \times d_x}$.
Applying Equations $\eqref{eq:admm_qp_z_simp}$ and $\eqref{eq:admm_qp_x_simp}$ allows us to express the ADMM iterations in a simplified form:

\begin{subequations} \label{eq:admm_qp_iter_simp}
\begin{align}
\begin{bmatrix}
\bx^{k+1}\\
\betta^{k+1}
\end{bmatrix}
& = -
\begin{bmatrix}
\bQ + \rho \bI_{\bx} &   \bA^T \\
\bA & 0
\end{bmatrix}^{-1}
\begin{bmatrix}
 \bp - \rho (\bz^{k} - \bmu^{k})\\
 -\blb
\end{bmatrix} \label{eq:admm_qp_x_iter_simp}\\
\bz^{k+1} &  = \Pi( \bx^{k+1}  + {\bmu}^k ) \label{eq:admm_qp_z_iter_simp}\\
{\bmu}^{k+1} & = {\bmu}^{k} +  \bx^{k+1} - \bz^{k+1} \label{eq:admm_qp_mu_iter_simp}
\end{align}
\end{subequations}

Observe that the per-iteration cost of the ADMM algorithm is dominated by solving the system of Equations $\eqref{eq:admm_qp_x_simp}$. Note, however, that this linear system is in general smaller than the Newton system found in a standard primal-dual interior-point solvers by a factor of approximately $5$, and therefore remains tractable for medium and large scale problems. Furthermore, if $\rho$ is static then the left-hand side matrix in Equation $\eqref{eq:admm_qp_x_simp}$ remains unchanged at each iteration and therefore is factorized only once at the onset of the ADMM algorithm. Furthermore, as we demonstrate below, this matrix factorization can subsequently be cached and invoked during backpropogation to compute the relevant gradients. Lastly, we note that if the matrices $\bQ$ and $\bA$ remain unchanged at each epoch of gradient descent, then the left-hand-side matrix needs to only be factorized \textbf{once} during the entire training process.

\subsubsection{ADMM-layer: unrolled differentiation}
Note that the standard unrolled differentiation approach `unrolls' the iterations in Equation $\eqref{eq:admm_qp_iter_simp}$  by standard backpropogation \citep{Hinton1986} and requires that each operation in Equation $\eqref{eq:admm_qp_iter_simp}$ be differentiable. We refer to  \citet{Domke2012} and \citet{Diamond2018} for a more comprehensive overview of unrolled differentiation. Our unrolled differentiation is invoked by the standard auto-differentiation routine in the \textit{torch} library.  

\subsubsection{ADMM-layer: KKT implicit differentiation}
As an alternative to unrolled differentiation, we note that the system of equations provided by the KKT conditions  at optimality is a fixed point mapping.  As outlined by \citep{Amos2017, Amos2019}, it is therefore possible to apply the implicit function theorem and derive the gradient of the primal-dual variables with respect to the problem variables in Program $\eqref{eq:pqp_main}$. In this section  we derive the KKT implicit differentiation for our ADMM solver. We begin with a few definitions.

\begin{defn}
Let $F \colon \mathbb{R}^{d_v} \times \mathbb{R}^{d_\theta} \to \mathbb{R}^{d_v}$ be a continuously differentiable function with variable $\bv$ and parameter $\btheta$. We define $\bv^*$ as a {\bf{fixed-point}} of $F$ at $(\bv^*,\btheta)$ if:
$$
F(\bv^*,\btheta) = \bv^*.
$$
\end{defn}

\begin{defn}
The {\bf{residual map}}, $G \colon \mathbb{R}^{d_v} \times \mathbb{R}^{d_\theta} \to \mathbb{R}^{d_v}$ of  a fixed point, $(\bv^*,\btheta)$, of $F$ is given by:
$$
G(\bv^*,\btheta) = F(\bv^*,\btheta) - \bv^* = 0.
$$
\end{defn}
The {\textbf{implicit function theorem}}, as defined by \citet{Dontchev2009},  then  provides the conditions on $G$ for which the Jacobian of the solution mapping with respect to $\btheta$ is well defined.

\begin{theorem}\label{th:if}
Let $G \colon \mathbb{R}^{d_v} \times \mathbb{R}^{d_\theta} \to \mathbb{R}^{d_v}$ be a continuously differentiable function in a neighborhood of $(\bv^*, \btheta)$ such that $G(\bv^*, \btheta) = 0$. Denote the nonsingular partial Jacobian of $G$ with respect to $\bv$ as $\nabla_{\bv} G(\bv^*, \btheta)$. Then $\bv(\btheta)$ is an implicit function of $\btheta$ and is continuously differentiable in a neighborhood, $\Theta$, of $\btheta$ with Jacobian:
\begin{equation} \label{eq:implicit}
\nabla_{\btheta} \bv(\btheta) = -[\nabla_{\bv} G(\bv(\btheta),\btheta) ]^{-1} \nabla_{\btheta} G(\bv(\btheta),\btheta) \quad \forall \quad \btheta \in \Theta.
\end{equation}

\end{theorem}

\begin{cor}\label{cor:if}
Let $F \colon \mathbb{R}^{d_v} \times \mathbb{R}^{d_\theta} \to \mathbb{R}^{d_v}$ be a continuously differentiable function with fixed-point $(\bv^*,\btheta)$. Then $\bv(\btheta)$ is an implicit function of $\btheta$ and is continuously differentiable in a neighborhood, $\Theta$, of $\btheta$ with Jacobian:
\begin{equation} \label{eq:implicit}
\nabla_{\btheta} \bv(\btheta) = [\bI_{\bv} - \nabla_{\bv} F(\bv(\btheta),\btheta) ]^{-1} \nabla_{\btheta} F(\bv(\btheta),\btheta) \quad \forall \quad \btheta \in \Theta.
\end{equation}
\end{cor}

For constrained quadratic programming, let us denote the primal-dual solution at optimality by $\bnu^* = (\bz^*, \tilde{\blambda}^*, \betta^*)$, where $\bz^*$ and $\betta^*$ are defined by Equations $\eqref{eq:admm_qp_x_simp}$ at optimality. Note that the dual variables associated with the inequality constraints are given by $\tilde{\blambda}^* = (\blambda^*_-,\blambda^*_+)$ with:
\begin{equation} \label{eq:lambda_admm}
\blambda^*_- = -\min(\rho \bmu^*,0) \quad \text{and} \quad  \blambda^*_+ = \max(\rho \bmu^*,0).
\end{equation}
We define the box inequality constraints as $\bG\bz \leq \bh$, where
$$
\bG = \begin{bmatrix} -\bI_{\bx} \\ \bI_{\bx} \end{bmatrix} \quad \text{and} \quad  \bh = \begin{bmatrix}-\bl \\ \bu \end{bmatrix}.
$$
Note that  all constraints are affine and therefore Slater's condition reduces to feasibility. The KKT conditions for stationarity, primal feasibility, and complementary slackness therefore defines a fixed-point at optimality $\bnu^*$ given by:

\begin{equation} \label{eq:kkt_G_qp}
G(\bnu^*,\btheta) =
\begin{bmatrix}
\bp + \bQ \bz^*  + \bG^T \tilde{\blambda}^* + \bA^T \betta^* \\
\diag(\tilde{\blambda}^*) (\bG \bz^* - \bh)\\
\bA\bz^* - \blb
\end{bmatrix}
=
\begin{bmatrix}
\bzero \\
\bzero \\
\bzero
\end{bmatrix}
\end{equation}
Applying Theorem \ref{th:if}, we take the differential of these conditions to give the following system of equations:

\begin{equation} \label{eq:diff}
\begin{split}
\begin{bmatrix}
 \bQ &  \bG^T &   \bA^T \\
\diag(\tilde{\blambda}^*)\bG & \diag (\bG \bz^* - \bh ) & 0\\
\bA & 0 & 0
\end{bmatrix}
\begin{bmatrix}
{\bdz }\\
{\bdlambda }\\
{\bdnu}
\end{bmatrix}
= -
\begin{bmatrix}
{\bdQ} \bz^* + {\bdp } + {\bdG}^T \tilde{\blambda}^* + {\bdA}^T \betta* \\
\diag(\tilde{\blambda} ^*) {\bdG } \bz^* - \diag(\tilde{\blambda} ^*) {\bdh}\\
{\bdA} \bz^* - {\bdb}
\end{bmatrix}.
\end{split}
\end{equation}

Observe that the left side matrix gives the optimality conditions of the convex quadratic problem, which, when solving by interior-point methods, must be factorized in order to obtain the solution to the nominal program \citep{Boyd2004}. The right side gives the differentials of the relevant functions at the achieved solution with respect to any of the problem variables. In practice, however, we never explicitly form the right-side Jacobian matrix directly. Instead we follow the work of \citet{Amos2017} and compute the left matrix-vector product of the Jacobian with the previous backward-pass gradient, $\frac{\partial \ell }{\partial \bz^*}$, as outlined below:

\begin{equation} \label{eq:diff_sol}
\begin{split}
\begin{bmatrix}
\bar{ \bd }_{\bz} \\
\bar{ \bd }_{\blambda} \\
\bar{ \bd }_{\betta}
\end{bmatrix}
= -
\begin{bmatrix}
\bQ &  \bG^T\diag(\tilde{\blambda}^*) &   \bA^T \\
\bG & \diag (\bG \bz^* - \bh ) & 0\\
\bA & 0 & 0
\end{bmatrix} ^{-1}
\begin{bmatrix}
\big( \frac{\partial \ell }{\partial \bz^*} \big)^T \\
0\\
0
\end{bmatrix}.
\end{split}
\end{equation}

Equation $\eqref{eq:diff_sol}$ allows for efficient computation of the  gradients with respect to any of the relevant problem variables. This is particularly true when using interior-point methods as the required gradients are effectively obtained `for free' upon factorization of the left matrix when obtaining the solution, $\bz^*$,  in the forward-pass. In the ADMM algorithm, however, we never explicitly solve the KKT system of equations and therefore we must form and factorize the left side KKT matrix during the backward pass routine. Observe that the left-side matrix in Equation $\eqref{eq:diff_sol}$ is of dimension $3d_z + d_{eq }$ and therefore for large scale problems solving this system of equations can be computationally burdensome. Finally, for the reader's interest, we state the gradients for all problem variables and refer the reader to \citet{Amos2017} for their derivation.

\begin{equation}\label{eq:kkt_partials}
\begin{aligned}
\frac{\partial \ell   }{\partial \bQ} & = \frac{1}{2} \Big(\bar{ \bd }_{\bz}  \bz^{*T} + \bz^* \bar{ \bd }_{\bz}^T \Big) & \qquad \frac{\partial \ell   }{\partial \bp} & = \bar{ \bd }_{\bz} \\
\frac{\partial \ell   }{\partial \bA} & =  \bar{ \bd }_{\betta}  \bz^{*T} + \betta^* \bar{ \bd }_{\bz} ^T   & \qquad \frac{\partial \ell   }{\partial \blb} & = -\bar{ \bd }_{\betta}  \\
\frac{\partial \ell   }{\partial \bG} & = \diag(\blambda ^*)  \bar{ \bd }_{\blambda}  \bz^{*T} + \blambda^* \bar{ \bd }_{\bz} ^T    & \qquad  \frac{\partial \ell   }{\partial \bh} & = - \diag(\blambda ^*) \bar{ \bd }_{\blambda}
\end{aligned}
\end{equation}

\subsubsection{ADMM-layer: fixed-point implicit differentiation}
In this section we demonstrate that the ADMM iterations in Equation $\eqref{eq:admm_qp_x_simp}$  can be cast as a fixed-point iteration of dimension $d_z + d_{eq}$. In many applications, $d_{eq}$ is typically much smaller than $d_z$, and therefore the proposed fixed-point differentiation will almost certainly decrease the computational overhead of the backward-pass routine. We begin with the following proposition. Note that  all proofs are available in the Appendix.

\begin{prop}\label{prop:admm_fixed_point}
Let $\bv^k = \bx^{k+1} +\bmu^k$, $\tilde{\bv}^k = (\bv^k,\betta^k)$ and define $F \colon \mathbb{R}^{d_{\tilde{v}}} \times \mathbb{R}^{d_\theta} \to \mathbb{R}^{d_{\tilde{v}}}$. Then the ADMM iterations in Equation $\eqref{eq:admm_qp_x_simp}$  can be cast as a fixed-point iteration of the form $F(\tilde{\bv},\btheta) = \tilde{\bv}$ given by:
\begin{align}\label{eq:admm_fp}
\begin{bmatrix}
\bv^{k+1}\\
\betta^{k+2}
\end{bmatrix}
& = -
\begin{bmatrix}
\bQ + \rho \bI_{\bv} &   \bA^T \\
\bA & 0
\end{bmatrix}^{-1}
\begin{bmatrix}
 \bp - \rho (2\Pi(\bv^{k}) - \bv^{k})\\
 -\blb
\end{bmatrix}
+
\begin{bmatrix}
\bv^{k}\\
\betta^{k+1}
\end{bmatrix}
-
\begin{bmatrix}
\Pi(\bv^{k})\\
\betta^{k+1}
\end{bmatrix}.
\end{align}
\end{prop}
We follow \citet{Boyd2018} and define the derivative of the projection operator, $\Pi$ as:
\begin{equation}\label{eq:project_box_derivative}
D\Pi(\bx) = \begin{cases}
                0 & \text{if } \bx_j < \bl_j\\
                1 & \text{if } \bl_j \leq \bx_j \leq \bl_j\\
                0 & \text{if } \bx_j > \bu_j\\
                \end{cases}.
\end{equation}
Observe that $D\Pi(\bx)$ is not continuously differentiable when $\bx_j = \bl_j$ or $\bx_j = \bu_j$. In practice, we can overcome the non-differentiability of $\Pi$ by introducing a small perturbation to $\bx$, thus moving $\bx$ away from the boundaries. Alternatively, smooth sigmoid based approximations to $\Pi(\bx)$ may also be suitable. In all experiments below, however, we invoke $D\Pi(\bx)$ directly as defined in Equation  $\eqref{eq:project_box_derivative}$.

The Jacobian, $\nabla_{\tilde{\bv}} F$, is therefore defined as:
\begin{align}\label{eq:jacob_F}
\nabla_{\tilde{\bv}} F
& = -
\begin{bmatrix}
\bQ + \rho \bI_{\bv} &   \bA^T \\
\bA & \bzero
\end{bmatrix}^{-1}
\begin{bmatrix}
- \rho (2D\Pi(\bv) - \bI_{\bv}) & 0\\
0 & 0
\end{bmatrix}
+
\begin{bmatrix}
\bI_{\bv} & 0\\
0 & \bI_{\betta}
\end{bmatrix}
-
\begin{bmatrix}
D\Pi(\bv) & 0\\
0 & \bI_{\betta}
\end{bmatrix}.
\end{align}
Corollary \ref{cor:if} therefore gives the desired Jacobian, $\nabla_{\btheta} \tilde{\bv}(\btheta)$, with respect to the parameter $\btheta$:
\begin{equation}\label{eq:jacob_v}
\nabla_{\btheta} \tilde{\bv}(\btheta) = [\bI_{\tilde{\bv}} - \nabla_{\tilde{\bv}} F(\tilde{\bv}(\btheta),\btheta) ]^{-1} \nabla_{\btheta} F(\tilde{\bv}(\btheta),\btheta)
\end{equation}
From the definition of $\bv$ we have that the Jacobians $\nabla_{\btheta} \bx(\btheta)$ and $\nabla_{\btheta} \betta(\btheta)$ are given by:
\begin{equation}\label{eq:jacob_v_etta}
\begin{bmatrix}
\nabla_{\btheta} \bx(\btheta)\\
\nabla_{\btheta} \betta(\btheta)
\end{bmatrix}
 =
\begin{bmatrix}
D\Pi(\bv) & 0\\
0 & \bI_{\betta}
\end{bmatrix}
\Big[\bI_{\tilde{\bv}} - \nabla_{\tilde{\bv}} F(\tilde{\bv}(\btheta),\btheta) \Big]^{-1} \nabla_{\btheta} F(\tilde{\bv}(\btheta),\btheta)
\end{equation}

As before we never form the Jacobians $\nabla_{\btheta} \bx(\btheta)$ and $\nabla_{\btheta} \betta(\btheta)$ directly. Instead, we compute the  left  matrix-vector  product  of  the  Jacobian  with  the  previous  backward-pass  gradient, $\frac{\partial \ell }{\partial \bz^*}$, as outlined below.

\begin{prop}\label{prop:admm_grads}
Let $\hat{ \bd }_{\bx}$ and $\hat{ \bd }_{\betta}$ be defined as:
\begin{equation}\label{eq:grads_admm}
\begin{split}
\begin{bmatrix}
\hat{ \bd }_{\bx}  \\
\hat{ \bd }_{\betta}
\end{bmatrix}
& =
\begin{bmatrix}
\bQ + \rho \bI_{\bv} &   \bA^T \\
\bA & 0
\end{bmatrix}^{-1}
\Big[\bI_{\tilde{\bv}} - \nabla_{\tilde{\bv}} F(\tilde{\bv}(\btheta),\btheta) \Big]^{-T}
\begin{bmatrix}
D\Pi(\bv) & 0\\
0 & \bI_{\betta}
\end{bmatrix}
\begin{bmatrix}
\big( - \frac{\partial \ell }{\partial \bz^*} \big)^T \\
0
\end{bmatrix}\\
& =
 \Bigg[
\begin{bmatrix}
D\Pi(\bv) & 0\\
0 & \bI_{\betta}
\end{bmatrix}
\begin{bmatrix}
\bQ + \rho \bI_{\bv} &   \bA^T \\
\bA & 0
\end{bmatrix}
+
\begin{bmatrix}
- \rho (2D\Pi(\bv) - \bI_{\bv}) & 0\\
0 & 0
\end{bmatrix}
\Bigg]^{-1}
\begin{bmatrix}
D\Pi(\bv) & 0\\
0 & \bI_{\betta}
\end{bmatrix}
\begin{bmatrix}
\big( - \frac{\partial \ell }{\partial \bz^*} \big)^T \\
0
\end{bmatrix}.
\end{split}
\end{equation}
Then the gradients of the loss function, $\ell$, with respect to problem variables $\bQ$, $\bp$, $\bA$ and $\blb$ are given by: .
\begin{equation}\label{eq:admm_partials}
\begin{aligned}
\frac{\partial \ell   }{\partial \bQ} & = \frac{1}{2} \Big(\hat{ \bd }_{\bx}   \bx^{*T} + \bx^* \hat{ \bd }_{\bx}^T \Big) & \qquad \frac{\partial \ell   }{\partial \bp} & = \hat{ \bd }_{\bx}  \\
\frac{\partial \ell   }{\partial \bA} & =  \hat{ \bd }_{\betta}  \bx^{*T} + \betta^* \hat{ \bd }_{\bx} ^T   & \qquad \frac{\partial \ell   }{\partial \blb} & = -\hat{ \bd }_{\betta}
\end{aligned}
\end{equation}
\end{prop}
Computing the gradients of the loss with respect to the box constraint variables, $\bl$ and $\bu$ is also straightforward.
\begin{prop}\label{prop:admm_grads_box}
Define $\tilde{\bmu}^*$ and $\hat{ \bd }_{\blambda}$ as:
\begin{equation}\label{eq:mu_tilde}
\tilde{\bmu}^*_j =  \begin{cases}
                \bmu^*_j & \text{if } \bmu^*_j \neq 0\\
                1 & \text{otherwise,} \\
                \end{cases}
\end{equation}
and
\begin{equation}\label{eq:d_lambda}
\hat{ \bd }_{\blambda} = \diag(\rho \tilde{\bmu}^*)^{-1} \Big( - \Big(\frac{\partial \ell }{\partial \bz^*} \Big)^T - \bQ \hat{ \bd }_{\bx} - \bA^T \hat{ \bd }_{\betta} \Big).
\end{equation}
Then the gradients of the loss function, $\ell$, with respect to problem variables $\bl$ and $\bu$ are given by:
\begin{equation}\label{eq:admm_partials_box}
\begin{aligned}
\frac{\partial \ell   }{\partial \bl} & = \diag(\blambda^*_-)\hat{ \bd }_{\blambda} & \qquad  \frac{\partial \ell   }{\partial \bu} & = -\diag(\blambda^*_+)\hat{ \bd }_{\blambda}.
\end{aligned}
\end{equation}
\end{prop}

We now have a framework for computing the gradient with respect to all problem variables by implicit differentiation of the fixed-point mapping of the transformed ADMM iterations. We re-iterate that the implicit differentiation of the KKT conditions requires solving a system of equations on the order of $3d_z + d_{eq}$. In contrast, the fixed-point iteration,  presented in Equation $\eqref{eq:admm_fp}$ is of dimension $d_z + d_{eq}$. As we will demonstrate shortly, reducing the dimension of the fixed-point mapping results in a considerable improvement in computational efficiency in the backward-pass.
\section{Computational experiments}\label{sec:results}
We present several experimental results that highlight the computational efficiency and performance accuracy of the ADMM-layer. In all experiments, computational efficiency is measured by the average runtime (in seconds), required to implement the forward-pass and backward-pass algorithms of each model. We compare across 4 models:
\begin{enumerate}
\item \textbf{ADMM Unroll:} ADMM in the forward-pass and unrolled differentiation in the backward-pass.
\item \textbf{ADMM KKT:} ADMM in the forward-pass and KKT implicit differentiation in the backward-pass.
\item \textbf{ADMM FP:} ADMM in the forward-pass and fixed-point implicit differentiation in the backward-pass.
\item \textbf{OptNet:} primal-dual interior-point method in the forward-pass and efficient KKT implicit differentiation in the backward-pass.
\end{enumerate}
Both the ADMM and interior-point solvers terminate when the $L_2$ norms of the primal and dual residuals are sufficiently small (i.e. less than some user-defined tolerance $\epsilon_{\text{tol}}$). In many applications, however, it is not always necessary to solve the batch QPs exactly during training. We therefore consider and compare the computational efficiency of each model over several stopping tolerances: $\epsilon_{\text{tol}} \in \{ 10^{-1}, 10^{-3}, 10^{-5}\}$. Going forward, the model label `ADMM FP 3', for example, denotes the ADMM FP model with a stopping tolerance of $10^{-3}$.

Lastly, first-order methods are known to be vulnerable to ill-conditioned problems and the resulting convergence rates can vary significantly when the data and algorithm parameters $(\rho)$ are poorly scaled. Many first-order solvers therefore implement a preconditioning and problem scaling initialization step (see for example \citep{Boyd2016, Schu2020, Stellato2020}). In our case, however, the QP problem variables are parameterized and are therefore expected to change at each epoch. Scaling and conditioning the problem variables at each epoch would potentially result in excessive computational overhead. As a result, our ADMM-layer implementation does not include a preconditioning step. Instead, in all experiments presented below, we normalize the problem data to have approximately unit standard deviation (on average) and manually scale the problem variables:  $\bp$, $\bQ$, $\bA$ and $\blb$ where appropriate. We find that for unit-scaled problem data, a value of $\rho \in \{0.10, 1.0\}$ provides a consistent rate of convergence. Indeed, an efficient and dynamic preconditioning and scaling routine is an interesting area of future research.

All experiments are conducted on an Apple Mac Pro computer (2.7 GHz 12-Core Intel Xeon E5,128 GB 1066 MHz DDR3 RAM) running macOS `Catalina'. All computations are run on an unloaded, single-threaded CPU. The software was written using the R programming language (version 4.0.0) and torch (version 0.6.0).

\subsection{Experiment 1: ADMM-layer performance}\label{sec:results_runtime}
 We conduct an experiment comparing the computational efficiency of the ADMM and OptNet models with various stopping tolerances. We randomly generate problem data of dimension: $$d_z \in \{10, 50, 100, 250, 500, 1000\}.$$ and for each trial implement the forward and backward algorithms on a mini-batch size of 128. Problem variables are generated as follows. We set $\bQ =  \frac{1}{2d_z}\bU^T\bU$ where entries of $\bU \in \mathbb{R}^{2d_z \times d_z}$ are sampled from a standard normal distribution. We randomly generate $\bp$ by sampling from the standard normal distribution and randomly generate $\bl$ and $\bu$ by randomly sampling from the uniform distribution with domain $[-2,-1]$ and $[1,2]$, respectively . Finally we set  $\bA = \bone$ and $\blb = 1$.

\begin{figure}[H]
  \centering
  \begin{subfigure}[b]{0.45\linewidth}
    \includegraphics[width=\linewidth , trim={0mm 0cm 0cm 0cm},clip]{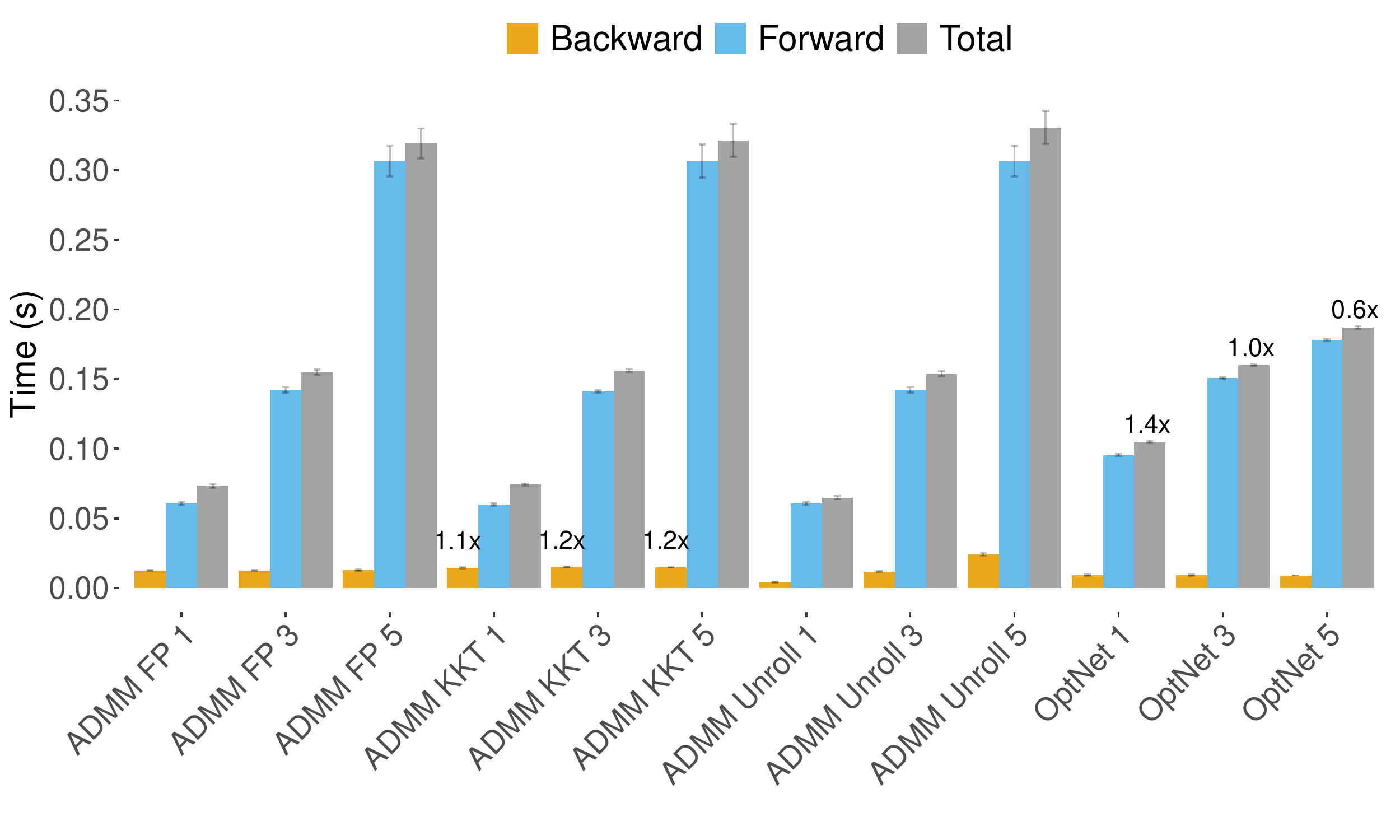}
    \caption{$d_z = 10$.}
  \end{subfigure}
  \begin{subfigure}[b]{0.45\linewidth}
    \includegraphics[width=\linewidth , trim={0mm 0cm 0cm 0cm},clip]{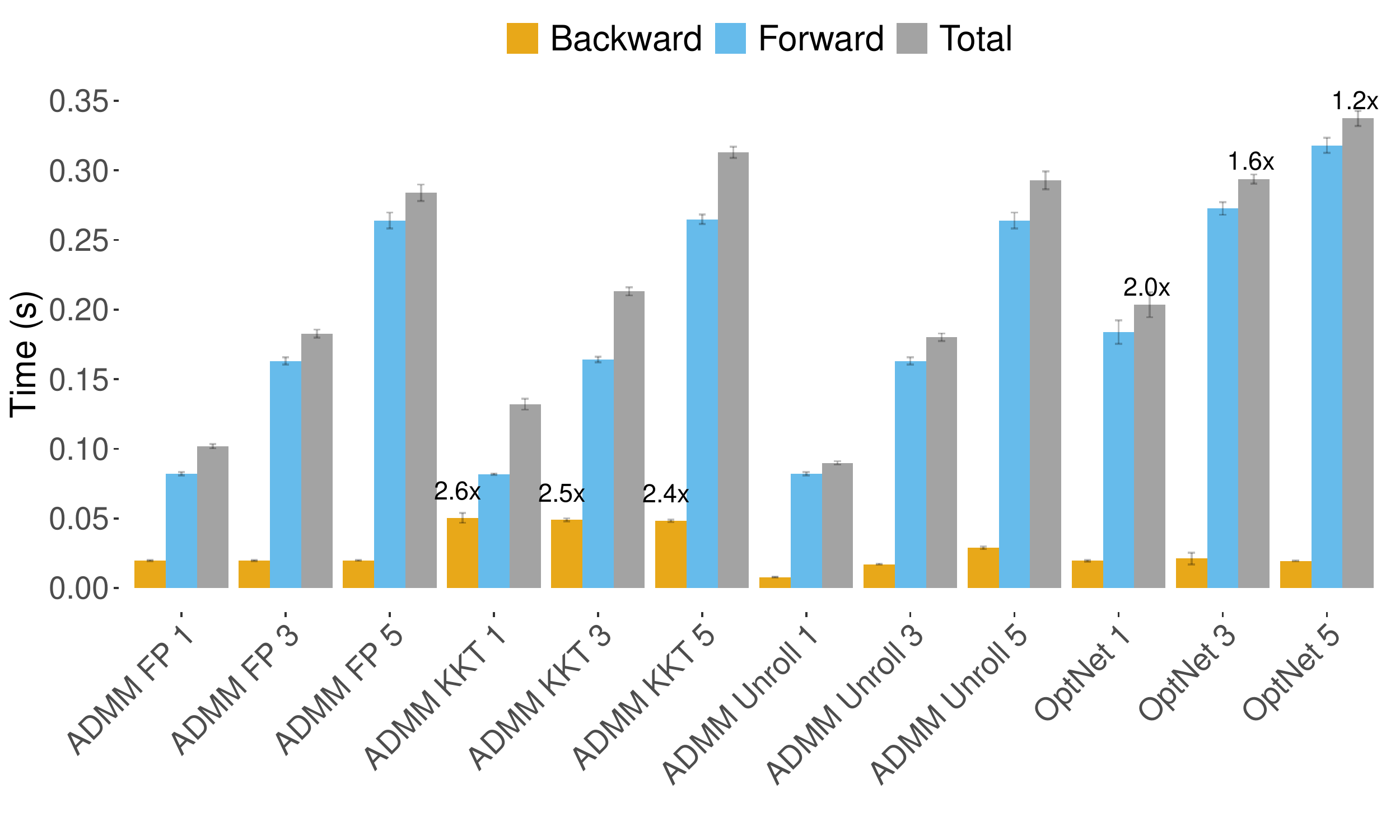}
    \caption{$d_z = 50$.}
  \end{subfigure}
    \begin{subfigure}[b]{0.45\linewidth}
    \includegraphics[width=\linewidth , trim={0mm 0cm 0cm 0cm},clip]{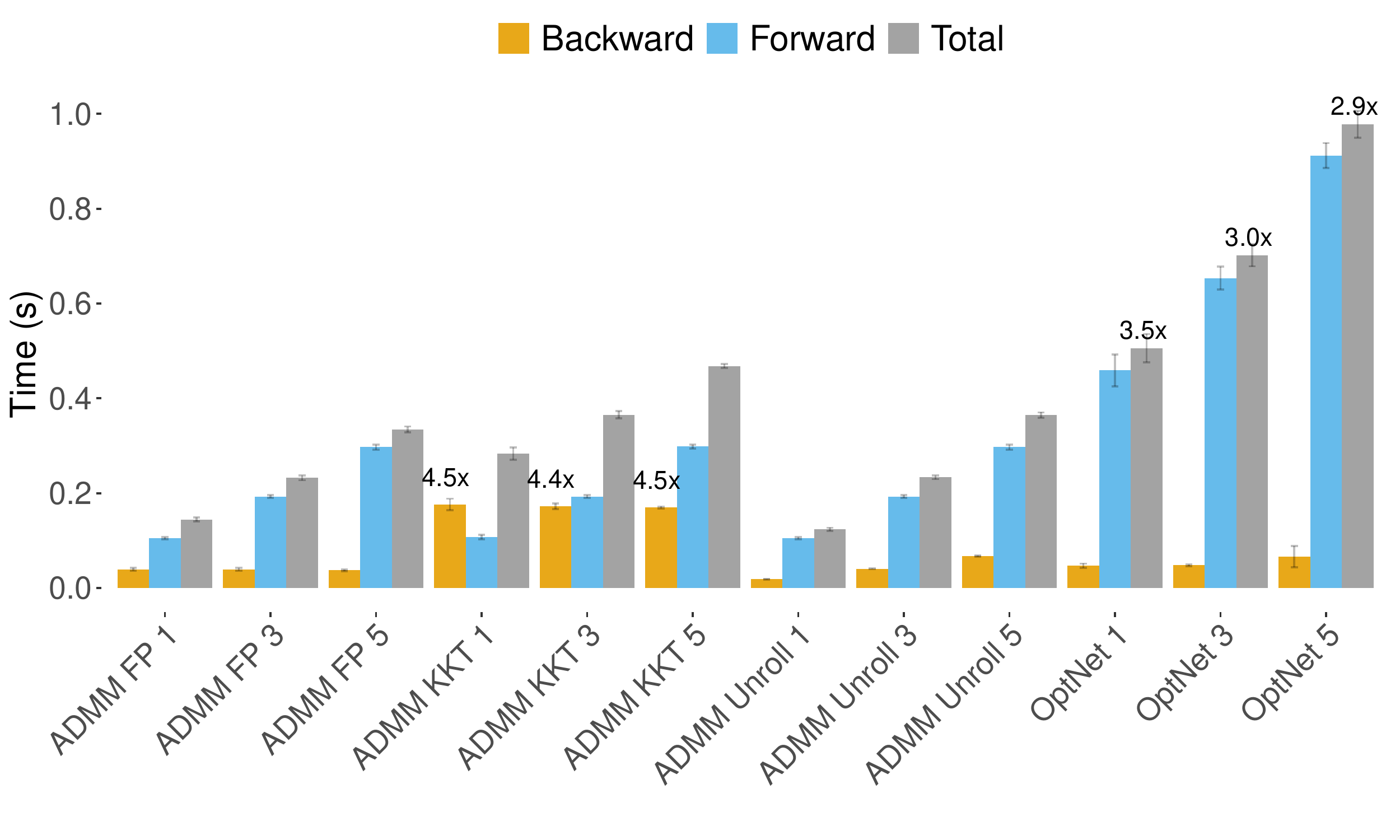}
    \caption{$d_z = 100$.}
  \end{subfigure}
  \begin{subfigure}[b]{0.45\linewidth}
    \includegraphics[width=\linewidth , trim={0mm 0cm 0cm 0cm},clip]{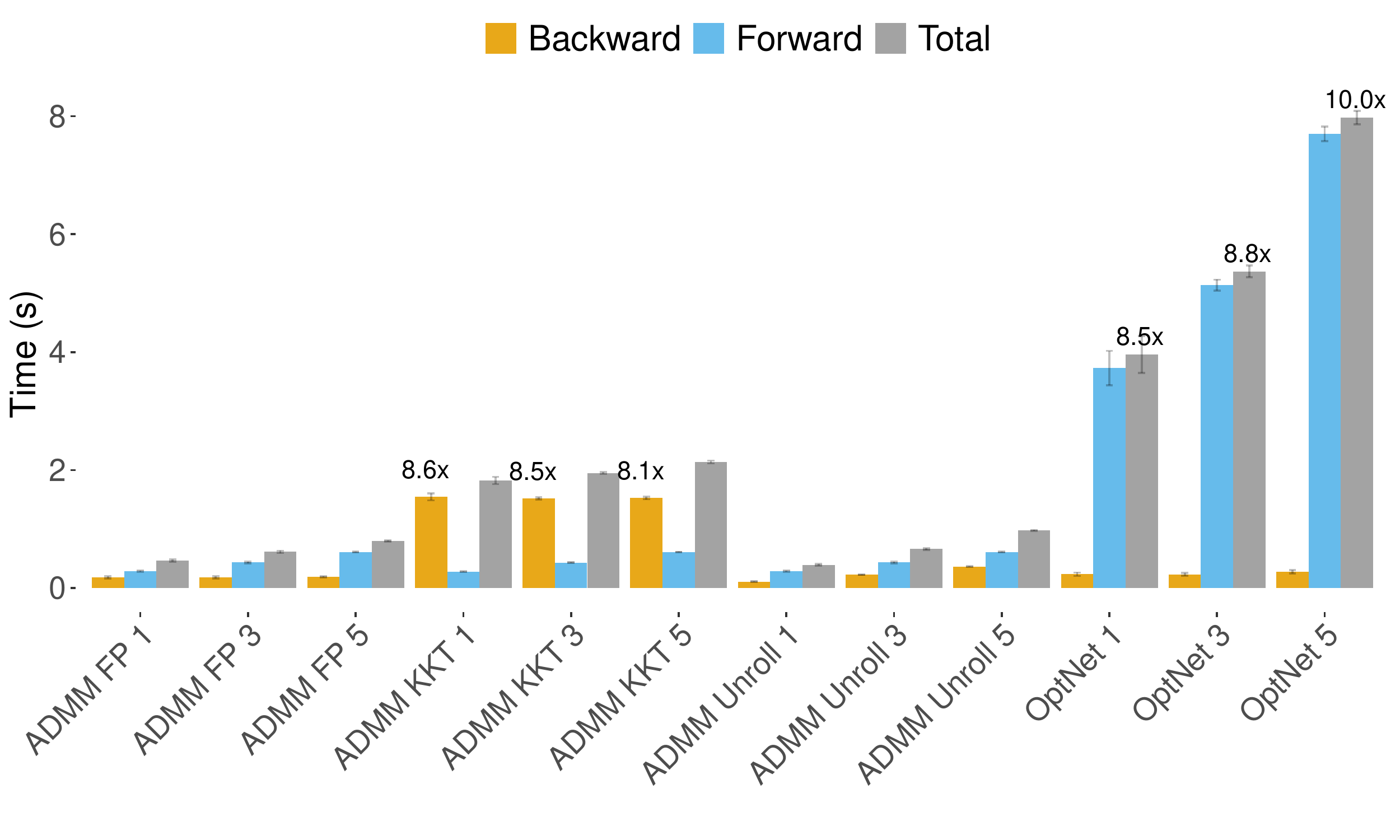}
    \caption{$d_z = 250$.}
  \end{subfigure}
  \begin{subfigure}[b]{0.45\linewidth}
    \includegraphics[width=\linewidth , trim={0mm 0cm 0cm 0cm},clip]{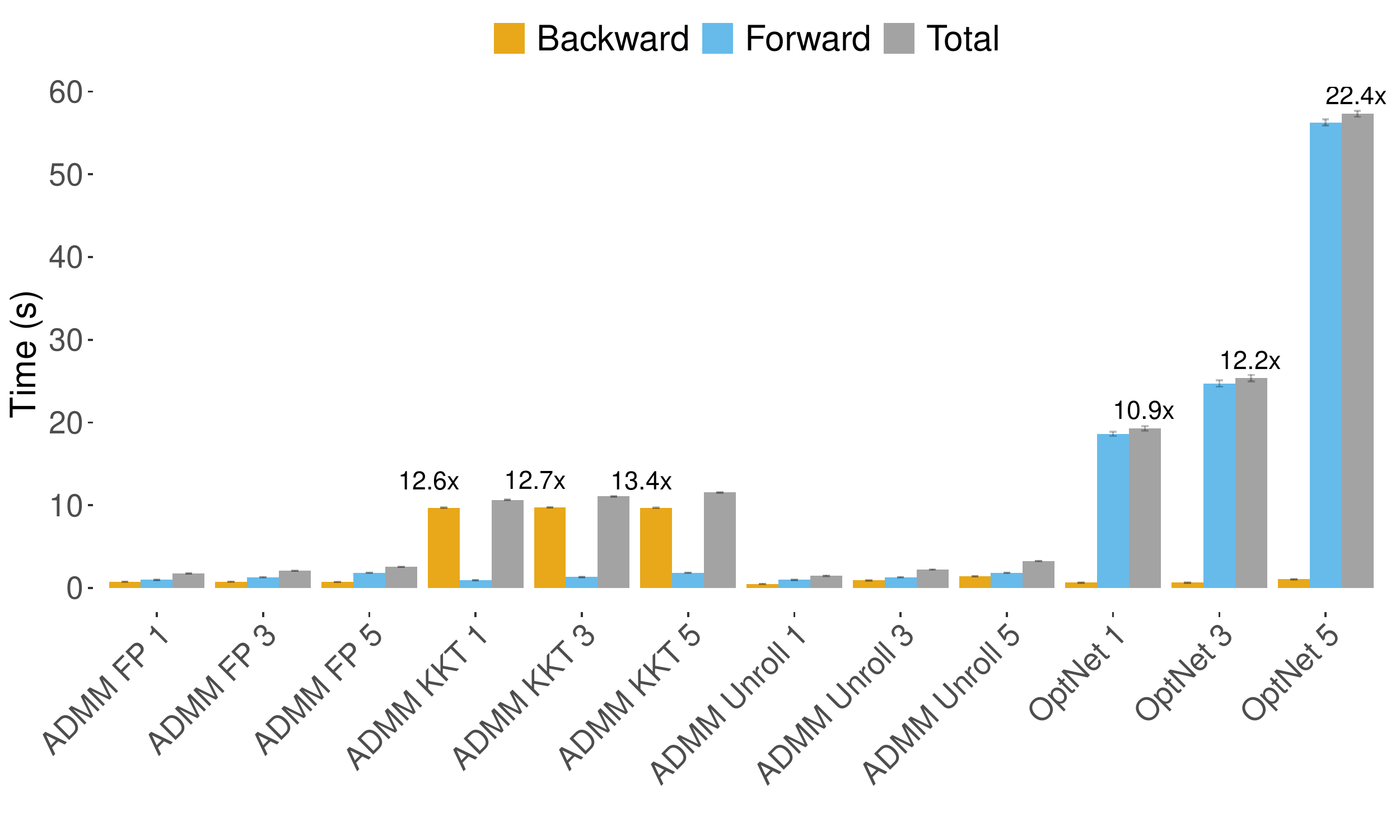}
    \caption{$d_z = 500$.}
  \end{subfigure}
  \begin{subfigure}[b]{0.45\linewidth}
    \includegraphics[width=\linewidth , trim={0mm 0cm 0cm 0cm},clip]{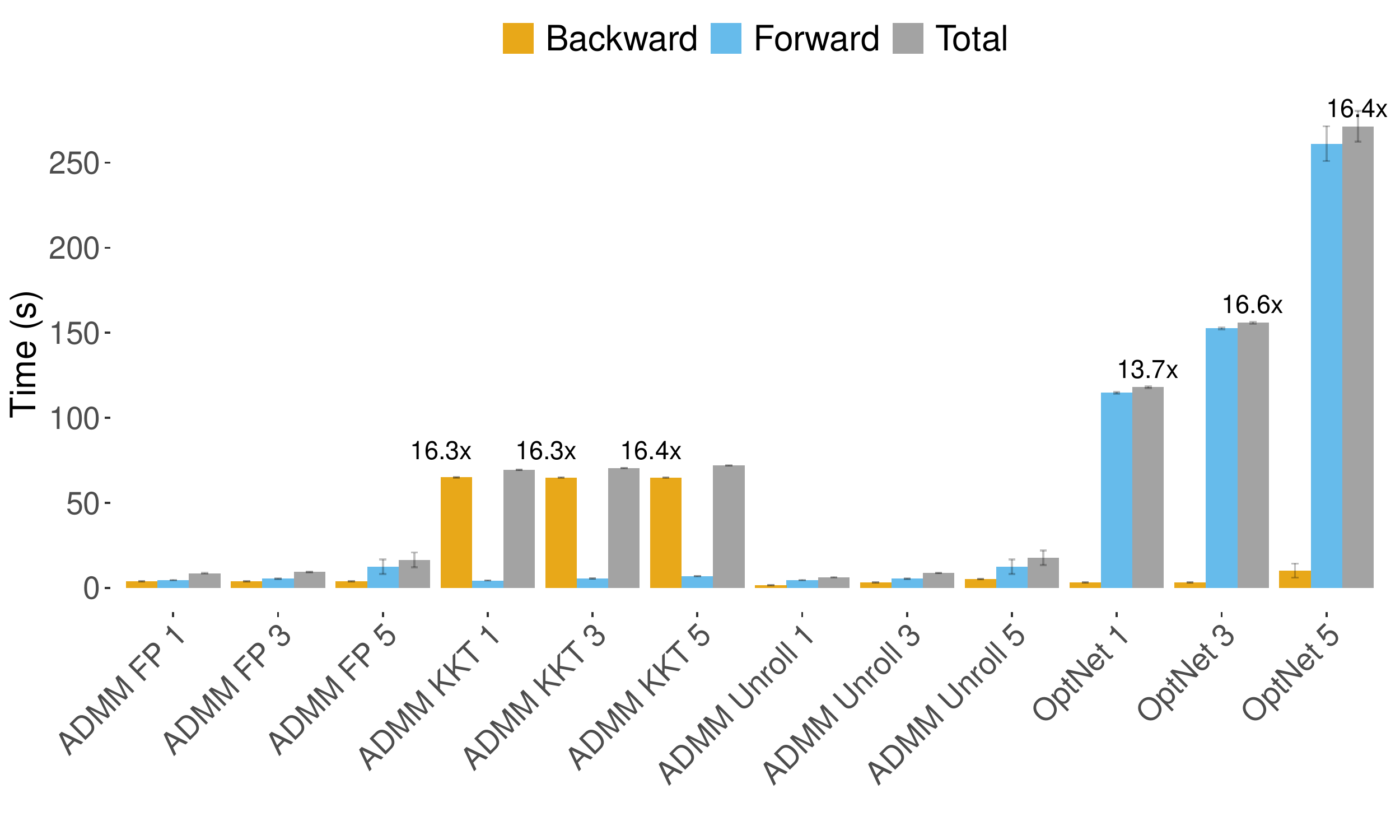}
    \caption{$d_z = 1000$.}
  \end{subfigure}
  \caption{Computational performance of ADMM-FP, ADMM-KKT, ADMM-Unroll and Optnet for various problem sizes, $d_z$, and stopping tolerances. Batch size $=128$.}
  \label{fig:admm_runtime}
\end{figure}

Figure \ref{fig:admm_runtime} provides the average runtime and $95\%$-ile confidence interval, evaluated over 10 trials, of the forward and backward-pass algorithms. We make several important observations. First, for small scale problems $(d_z < 100)$, there is negligible performance differences across all methods of the same stopping tolerance. As expected, the total runtime increases as the required stopping tolerance decreases. For medium scaled problems $(100 < d_z < 1000)$ we observe a substantial performance degradation in both the ADMM-KKT models and the OptNet models, in comparison to the ADMM-FP and ADMM-Unroll models. Specifically, the ADMM-KKT model exhibits an increase in computation time that is anywhere from $4$ to $16$ times larger than the corresponding ADMM fixed-point backward-pass implementation. This result is not surprising as the ADMM-KKT backward-pass algorithm must first form and then factorize the KKT system of equations, which is of dimension $3d_z + d_{eq}$. In contrast, the fixed-point backward-pass routine solves a system of equations of size $d_z + d_{eq}$ and is shown to be comparable in computational efficiency to the OptNet backward-pass algorithm. Furthermore, for problems of size $d_z \geq 250$, we note that the ADMM-FP and ADMM-Unroll models are approximately an order of magnitude more efficient than the corresponding OptNet models. For example, when $d_z = 1000$ and $\epsilon_{\text{tol}} = 10^{-3}$, the total runtime for the OptNet model is $150$ seconds, whereas the total runtime for the ADMM model is less than $10$ seconds; over an order of magnitude faster. This increase in computational performance will ultimately enable training architectures that can practically support substantially larger quadratic optimization problems. Lastly, we note that while the ADMM-unroll algorithm is relatively efficient, it requires a significantly larger memory footprint, which may be impractical in some settings. Furthermore, we observe that the as the stopping tolerance decreases the computation time of the unrolled backward-pass increases. In contrast, the fixed-point implicit differentiation method is invariant to the number of `inner' iterations performed in the ADMM forward-pass. Going forward, we choose to work with the ADMM fixed-point model as it is efficient from both a computational and memory standpoint.  


\subsection{Experiment 2: learning $\bp$}\label{sec:results_learn_p}
We now consider a full training experiment whereby the objective is to learn a parameterized model for the variable $\bp$, that is optimal in the context of the remaining QP problem variables. This problem was considered by \citet{Donti2017} with applications to power scheduling and battery storage, and more recently by  \citet{Butler2020IPO} for optimal return forecasting within the context of a mean-variance portfolio.  We refer to the aforementioned work for more details.

The learning process can be posed as a bi-level optimization program where  the objective is to learn a parameter $\btheta$ in order to minimize the average QP loss induced by the optimal decision policies $\{\toi{\bz^*(\btheta)}\}_{i=1}^m$. Program $\eqref{eq:qp_p}$ is referred to as an integrated predict and optimize (IPO) model as the prediction model for $\bp$ is fully integrated with the resulting down-stream decision-based optimization model.
\begin{equation} \label{eq:qp_p}
\begin{split}
\minimize_{ \btheta} \quad & \frac{1}{m} \sum_{i = 1}^m \Big(  \toit{\bz^*(\btheta)} \toi{\bp} + \frac{1}{2} \toit{\bz^*(\btheta)} \toi{\bQ} \toi{\bz^*(\btheta)}  \Big)  \\
\text{subject to }\quad  &     \toi{\bz^*(\btheta)}   = \argmin_{\bz}  -\bz^T \toi{\hat{\bp}(\btheta)}  + \frac{1}{2} \bz^T \toi{\hat{\bQ}} \bz \quad \forall i=1,...,m\\
\quad & \bA \toi{\bz^*(\btheta)} = \blb \quad \forall i=1,...,m\\
\quad &  \bl \leq \toi{\bz^*(\btheta)} \leq \bu \quad \forall i=1,...,m.
\end{split}
\end{equation}
Here $\toi{\bp}$ and $\toi{\bQ}$ denote the ground truth problem data and are generated as follows. We let ${\toi{\bp} \sim \mathcal{N}(\toit{\bw} \btheta_0 + \tau \toi{\bepsilon}, \bQ)}$ where $\bQ \in \mathbb{R}^{d_z \times d_z}$ has entry $(j, k)$ equal to $\rho_{\bp}^{|j-k|}$. We set $\rho_{\bp}= 0.50$ and generate the auxiliary feature data from the standard normal distribution, $\toi{\bw} \sim \mathcal{N}(\bm 0, \bI_{\bw})$.  The residuals, $\toi{\bepsilon} \sim \mathcal{N}(\bm 0, \bQ)$, preserve the desired correlation structure and the scalar value $\tau$ controls the signal-to-noise level. All experiments target a signal-to-noise level of $0.10$.

We let $\toi{\hat{\bp}(\btheta)}$ denote the estimate of $\toi{\bp}$ according to the linear model:
\begin{equation}
\toi{\hat{\bp}(\btheta)} = \toit{\bw} \btheta.
\end{equation}
The bound constraints, $\bl$ and $\bu$, are generated by randomly sampling from the uniform distribution with domain $[-1,0]$ and $[0,1]$, respectively and we set  $\bA = \bone$ and $\blb = 1$. In all experiments we set the stopping tolerance to $\epsilon_{\text{tol}} = 10^{-3}$.

\begin{figure}[H]
  \centering
  \begin{subfigure}[b]{0.45\linewidth}
    \includegraphics[width=\linewidth, trim={0cm 0cm 0cm 0cm},clip]{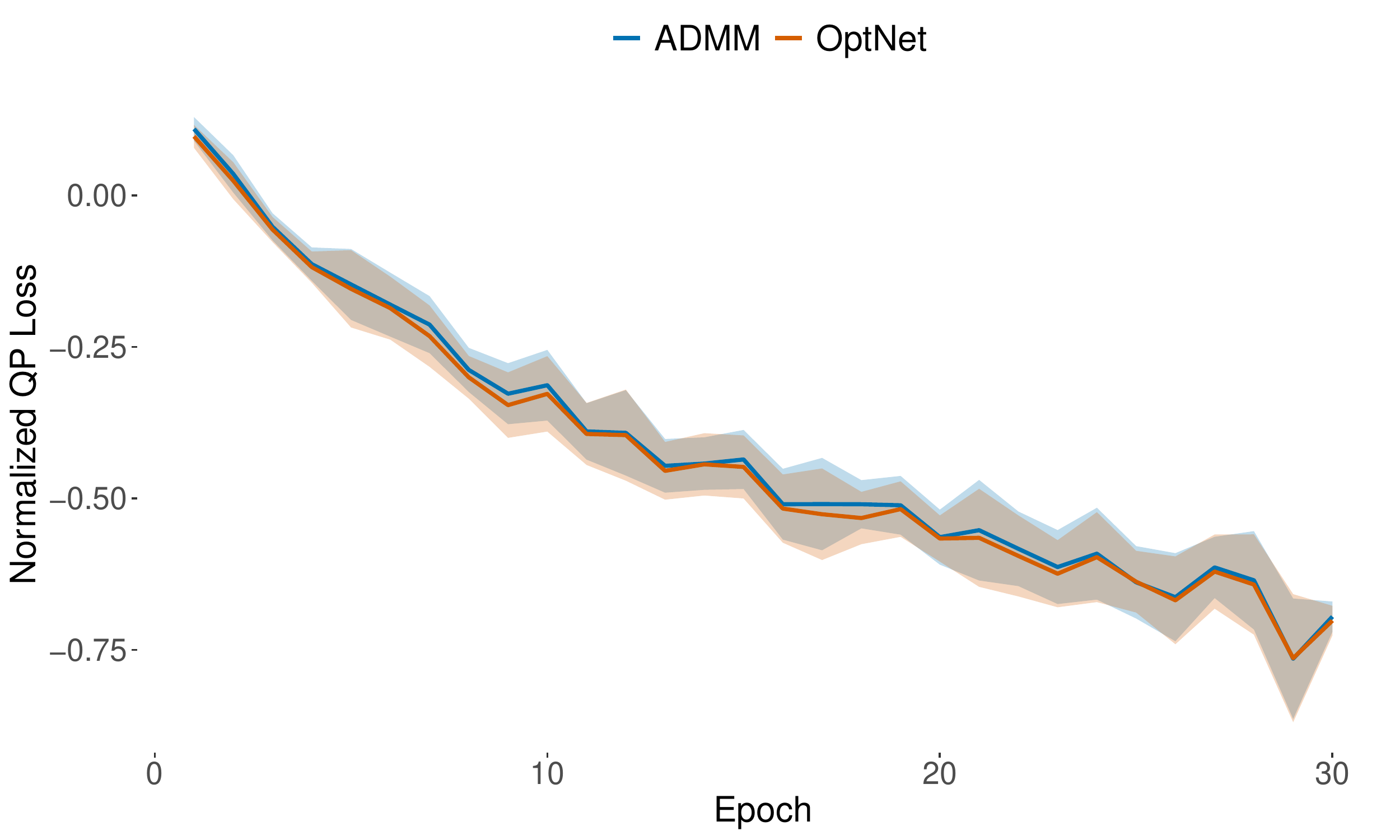}
    \caption{Training Loss.}
  \end{subfigure}
  \begin{subfigure}[b]{0.45\linewidth}
    \includegraphics[width=\linewidth, trim={0cm 0cm 0cm 0cm},clip]{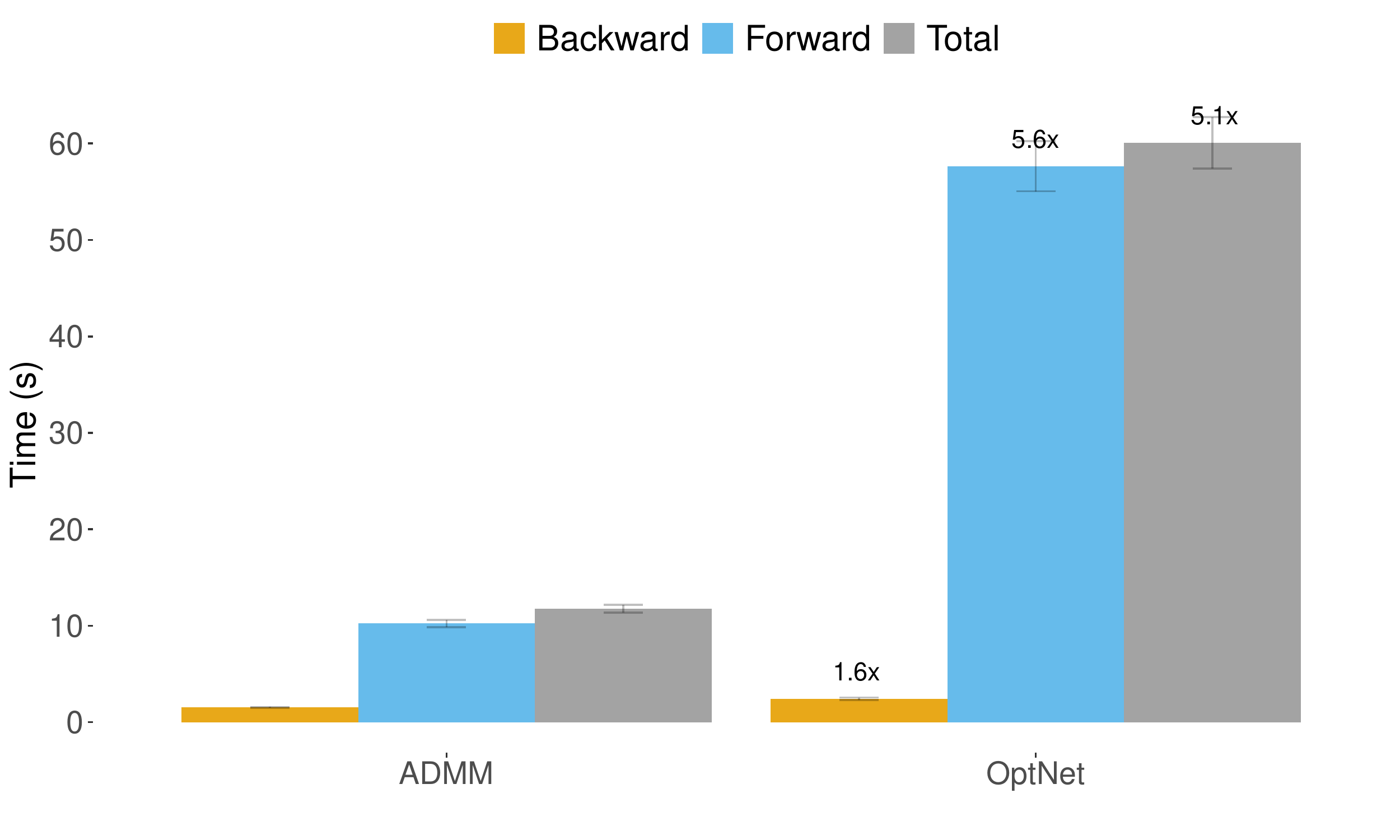}
    \caption{Computational Performance.}
  \end{subfigure}
  \caption{Training loss and computational performance for learning $\bp$. Batch size $=32$ and $d_z = 250$.}
  \label{fig:learn_p_250}
\end{figure}

\begin{figure}[H]
  \centering
  \begin{subfigure}[b]{0.45\linewidth}
    \includegraphics[width=\linewidth, trim={0cm 0cm 0cm 0cm},clip]{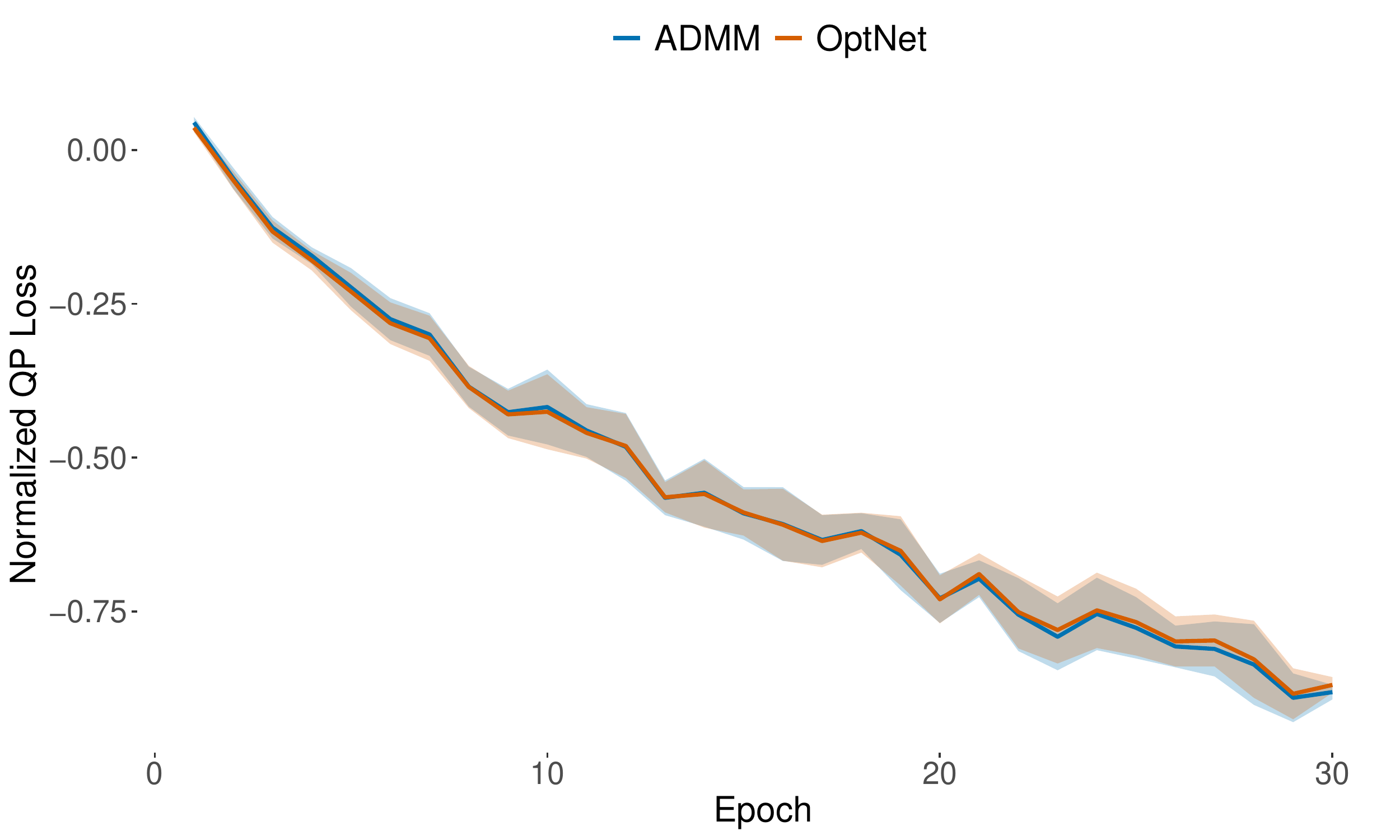}
    \caption{Training Loss}
  \end{subfigure}
  \begin{subfigure}[b]{0.45\linewidth}
    \includegraphics[width=\linewidth, trim={0cm 0cm 0cm 0cm},clip]{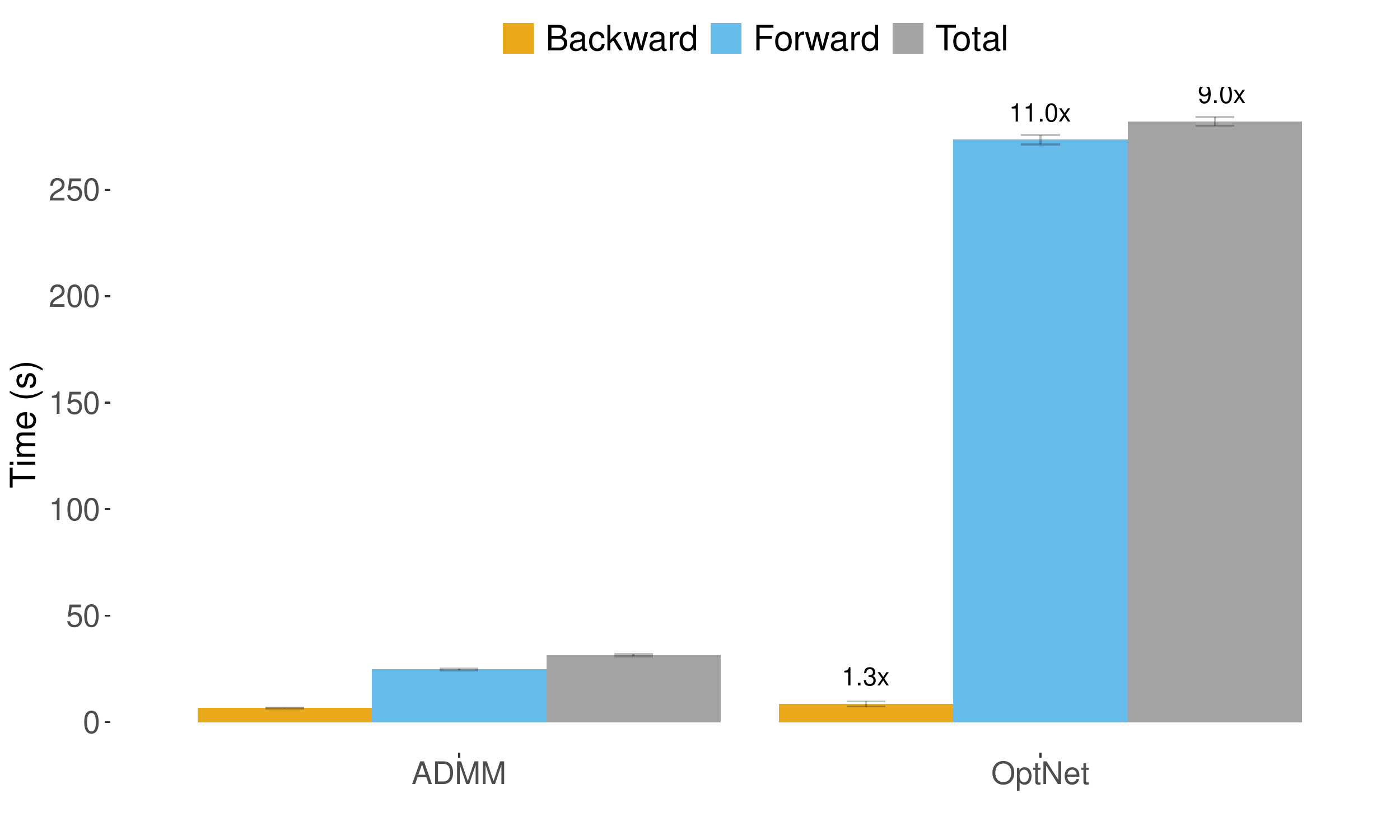}
    \caption{Computational Performance}
  \end{subfigure}
\caption{Training loss and computational performance for learning $\bp$. Batch size $=32$ and $d_z = 500$.}
  \label{fig:learn_p_500}
\end{figure}

\begin{figure}[H]
  \centering
  \begin{subfigure}[b]{0.45\linewidth}
    \includegraphics[width=\linewidth, trim={0cm 0cm 0cm 0cm},clip]{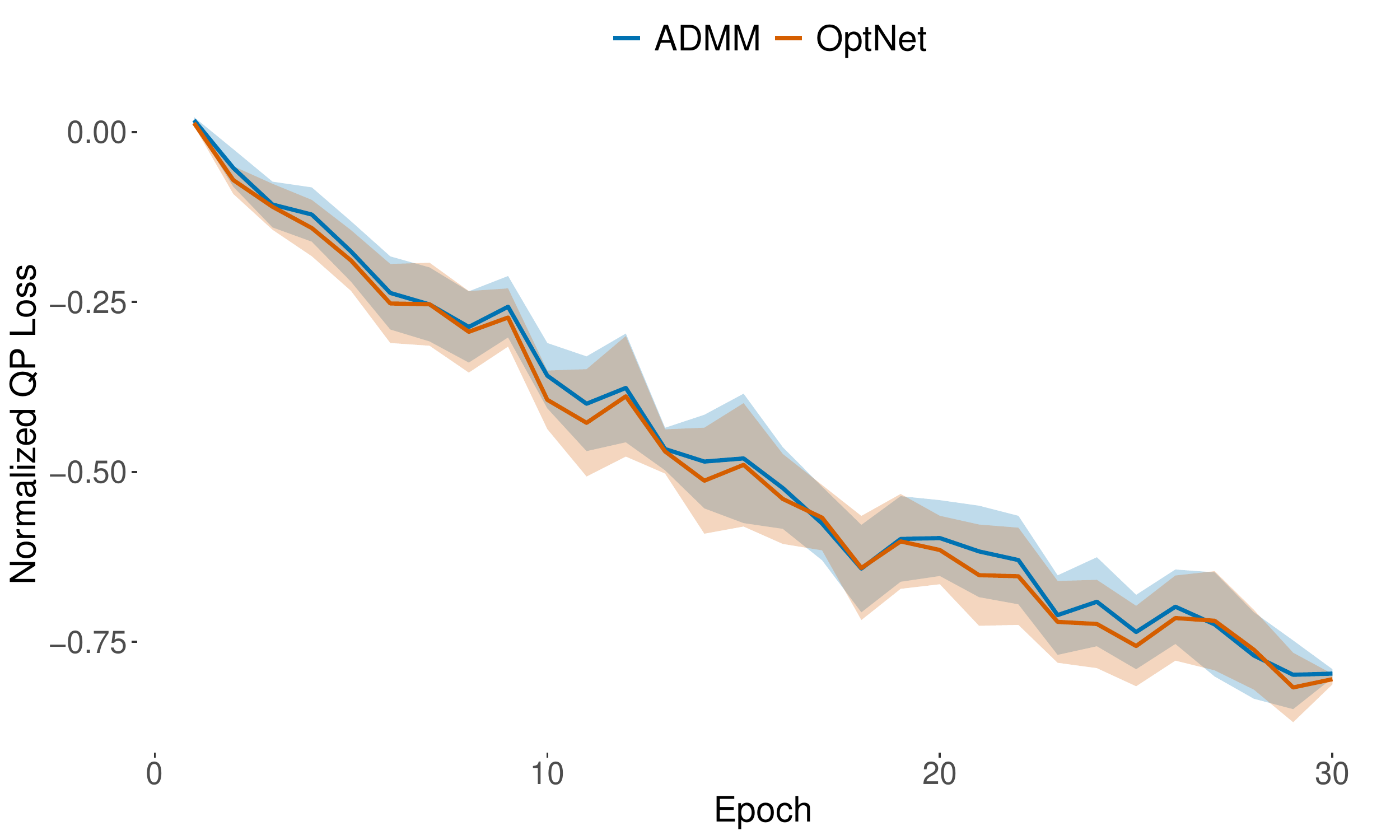}
    \caption{Training Loss}
  \end{subfigure}
  \begin{subfigure}[b]{0.45\linewidth}
    \includegraphics[width=\linewidth, trim={0cm 0cm 0cm 0cm},clip]{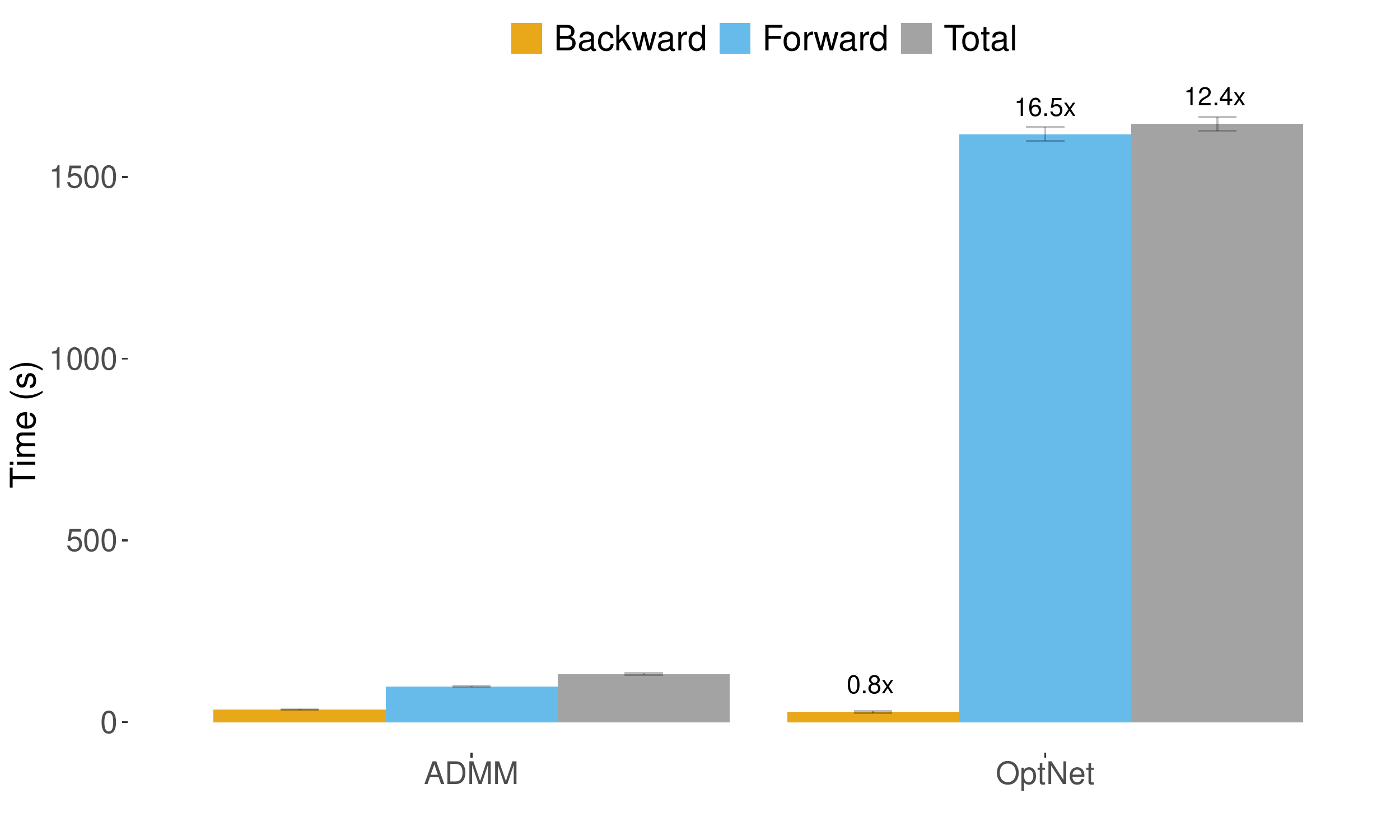}
    \caption{Computational Performance}
  \end{subfigure}
  \caption{Training loss and computational performance for learning $\bp$. Batch size $=32$ and $d_z = 1000$.}
  \label{fig:learn_p_1000}
\end{figure}

We randomly generate problem data of dimension $d_z \in \{250,500,1000\}$. The training process for each trial consists of $30$ epochs with a mini-batch size of $32$. Figures \ref{fig:learn_p_250}(a) - \ref{fig:learn_p_1000}(a) report the average training loss at each epoch and the $95\%$-ile confidence interval, evaluated over $10$ independent trials. Observe that the loss curves for the ADMM model and OptNet model are almost identical at each epoch, suggesting that the training accuracy of the ADMM model and OptNet model are equivalent. Conversely, Figures \ref{fig:learn_p_250}(b) - \ref{fig:learn_p_1000}(b), compare the average and $95\%$-ile confidence interval time spent executing the forward and backward pass algorithms. When $d_z = 250$ the ADMM model is shown to be approximately $5$ times faster than the OptNet model. Furthermore, when $d_z = 500$ and $d_z = 1000$, the ADMM model is a full order of magnitude faster than the OptNet model. More concretely, when $d_z = 1000$ the entire learning process takes approximately $1600$ seconds to train the OptNet model, but less than $130$ seconds to train the ADMM model to an equal level of accuracy.

\subsection{Experiment 3: learning $\bA$}\label{sec:results_learn_A}
We now present a real-world experiment from portfolio optimization whereby the objective is to learn a parameterized model for the variable $\bA$. We consider an asset universe of $d_z = 255$ liquid US stocks traded on major U.S. exchanges (NYSE, NASDAQ, AMEX, ARCA). The universe is summarized in Table \ref{table:stocks}, with representative stocks from each of the Global Industry Classification Standard (GICS) sectors. Weekly price data is given from January $1990$ through December $2020$, and is provided by Quandl.

We denote the matrix of weekly return observations as $\bA=[\ba^{(1)},\ba^{(2)},...,\ba^{(m)} ] \in \mathbb{R}^{m\times d_z}$ with $m > d_z$. Let $\toi{\bQ}  \in \mathbb{R}^{d_z \times d_z}$ denote the symmetric positive definite covariance matrix of asset returns. We define the portfolio $\toi{\bz} \in \mathbb{R}^{d_z}$,  where the element, $\toi{\bz}_j$, denotes the proportion of total capital invested in the $j^{\text{th}}$ asset at time $i$.

We define the Sharpe ratio at observation $i$ as the ratio of portfolio return to portfolio risk, where risk is measured by the portfolio volatility (standard deviation).
\begin{equation}
\toi{S_R }= \frac{ \toit{\ba} \toi{\bz} }{ \sqrt{ \toit{\bz} \toi{\bQ} \toi{\bz} }}
\end{equation}
We consider a long-only $(\toi{\bz} \geq 0)$, fully invested $(\bone^T\toi{\bz} = 1)$  max-Sharpe portfolio optimization, presented in Program $\eqref{eq:max_sharpe}$:

\begin{equation}\label{eq:max_sharpe}
\begin{split}
\maximize_{\bz} \quad & \frac{\toit{\ba}\bz }{ \sqrt{ \bz^T \toi{\bQ} \bz }}\\
\text{subject to }\quad  & \bone^T\bz = 1, \quad 0 \leq \bz \leq 1\\
\end{split}
\end{equation}
Observe, however, that the Sharpe ratio is not convex in $\bz$ but is homogeneous of degree zero. We follow \citet{Corn2006} and re-cast Program $\eqref{eq:max_sharpe}$ as a convex quadratic optimization program:
\begin{equation}\label{eq:max_sharpe_2}
\begin{split}
\minimize_{\bz} \frac{1}{2} \quad & \bz^T \toi{\bQ} \bz \\
\text{subject to }\quad  & \toit{\ba}\bz = 1, \quad   \bz \geq 0.\\
\end{split}
\end{equation}
Note that the fully-invested constraint can be enforced by normalizing the optimal weights, $\bz^*$.

As before, the learning process can be posed as a bi-level optimization program where  the objective is to learn a parameter $\btheta$ and the associated constraints, $\toi{\hat{\ba}(\btheta)}$, in order to maximize the average realized Sharpe ratio induced by the optimal decision policies $\{\toi{\bz^*(\btheta)}\}_{i=1}^m$.

\begin{equation} \label{eq:max_sharpe_bi}
 \begin{split}
 \minimize_{ \btheta } \quad & -\frac{1}{m} \sum_{i = 1}^m  \frac{\toit{\ba} \bz^*(\btheta) }{ \sqrt{ \toit{\bz^*(\btheta)} \toi{\bQ} \toi{\bz^*(\btheta)} } }    \\
 \text{subject to }\quad  &  \bz^*(\btheta)   = \argmin_{\bz}  \frac{1}{2} \bz^T \toi{\bQ} \bz \quad \forall i=1,...,m\\
 \quad &\toit{\hat{\ba}(\btheta)} \toi{\bz^*(\btheta)} =1, \quad \bone^T \toi{\bz^*(\btheta)} =1, \quad \toi{\bz^*(\btheta)}  \geq 0 \quad \forall i=1,...,m
 \end{split}
 \end{equation}

In reality, we do not know the true value of $\toi{ \ba }$ at decision time and instead we estimate $\toi{ \ba }$ through associated auxiliary feature variables $\toi {\bw } \in \mathbb{R}^{d_w}$. Again we consider a linear   model of the form:
\begin{equation}\label{eq:ff5_ret_forecast}
\toi {\hat{\ba} } =  \btheta^T \toi {\bw}.
\end{equation}
In this experiment, asset returns, $\toi{ \ba }$  are modelled using the well-known Famma-French Five (FF5) factor model \citep{Fama2015}, provided by the Kenneth R. French data library.

The goal is to observe the training and out-of-sample performance of the ADMM model in comparison to the OptNet model. As a benchmark, we include the out-of-sample performance of an equally weighted portfolio, and a max-Sharpe portfolio where $\btheta$ is fit by ordinary least-squares (OLS).  All experiments are trained on data from January $1990$ through December $2014$. The out-of-sample period begins in January $2015$ and ends in December $2020$. Portfolios are formed at the close of each week, and rebalanced on a weekly basis.

\begin{figure}[H]
  \centering
  \begin{subfigure}[b]{0.45\linewidth}
    \includegraphics[width=\linewidth, trim={0cm 0cm 0cm 0cm},clip]{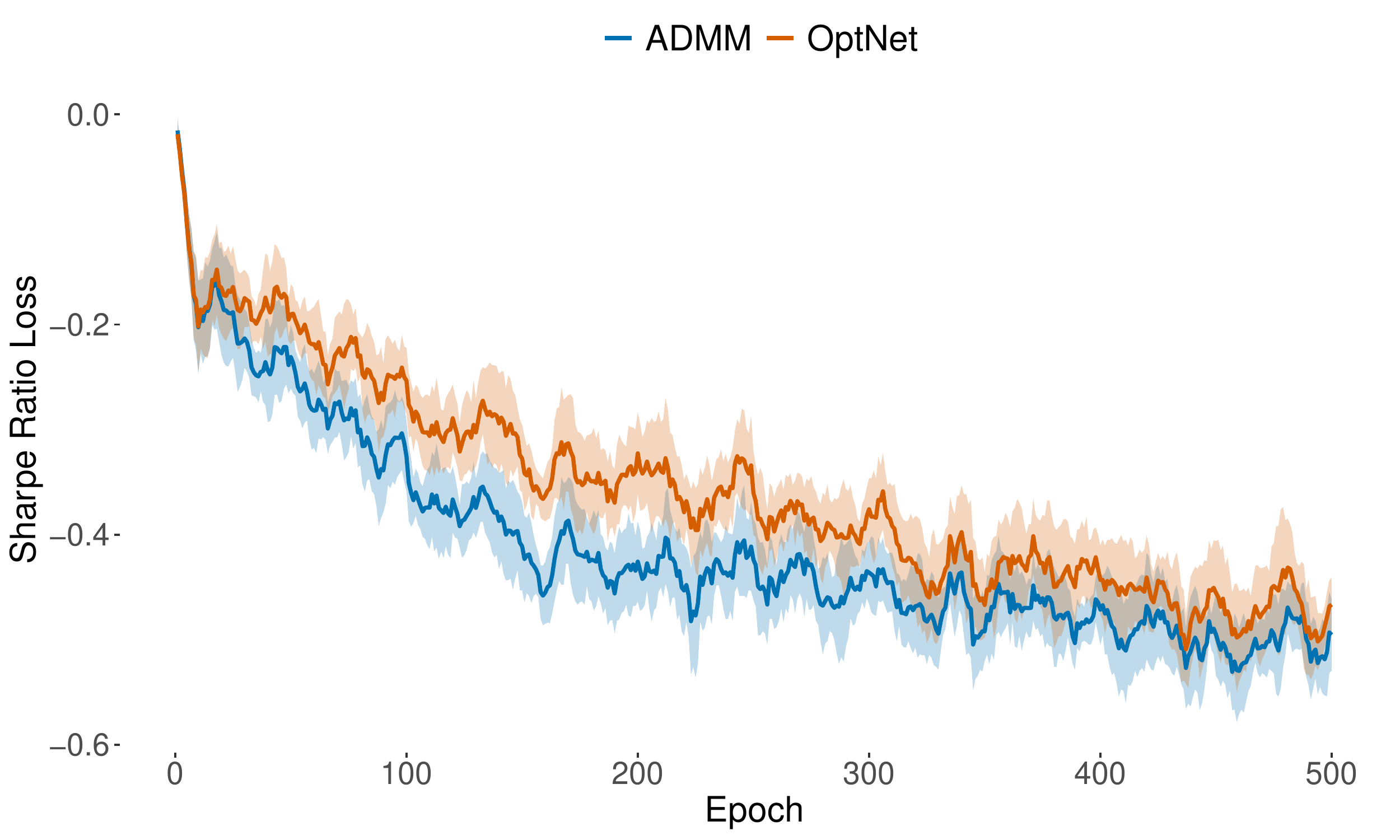}
    \caption{Training Loss}
  \end{subfigure}
  \begin{subfigure}[b]{0.45\linewidth}
    \includegraphics[width=\linewidth, trim={0cm 0cm 0cm 0cm},clip]{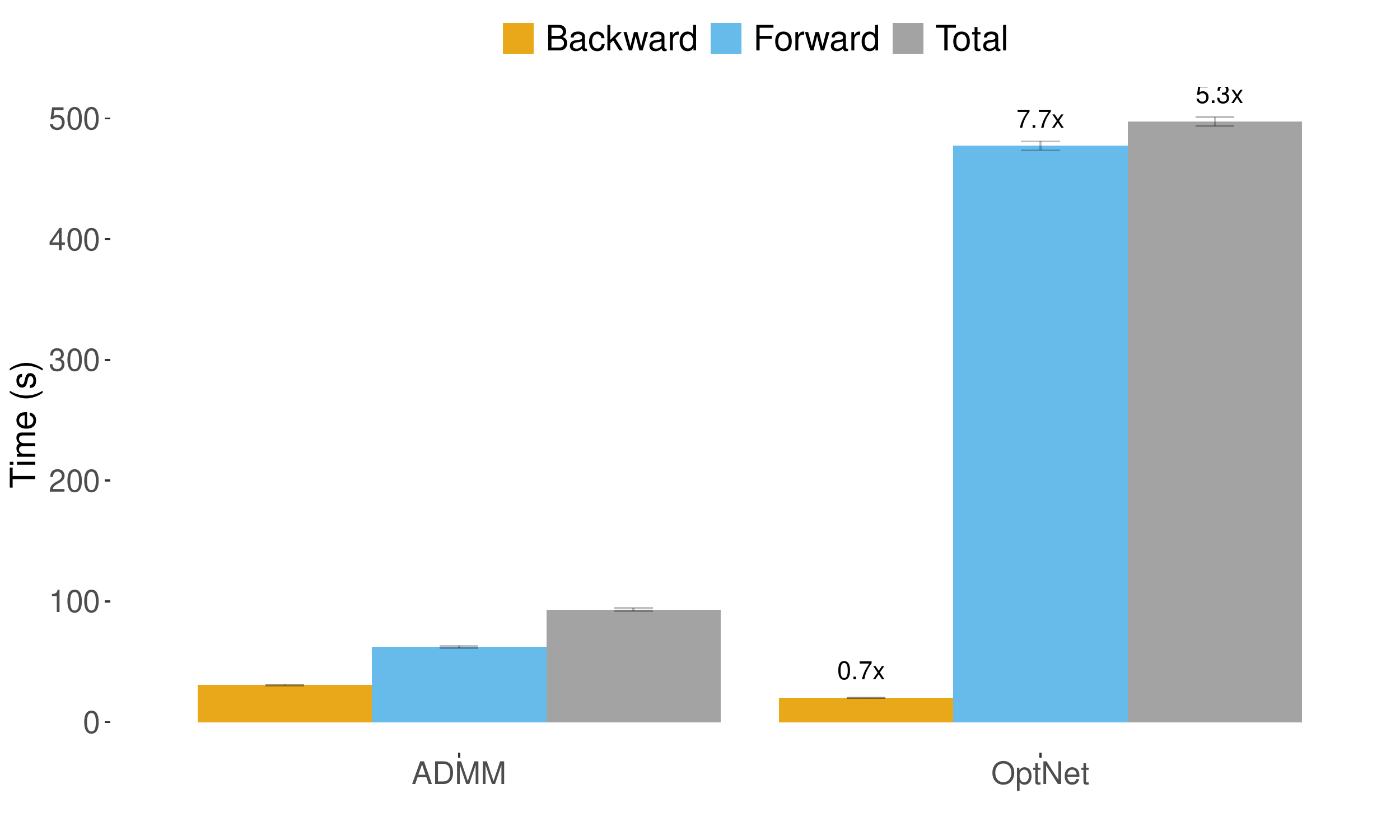}
    \caption{Computational Performance}
  \end{subfigure}
  \caption{Training loss and computational performance for learning $\bA$ on US stock data. \\ Batch size $=32$ and $d_z = 255$.}
  \label{fig:learn_A}
\end{figure}

 The training process for each trial consists of $500$ epochs with a mini-batch size of $32$. Portfolio models are fit to an accuracy of $\epsilon_{\text{tol}} = 10^{-4}$ in training, and a higher accuracy of $\epsilon_{\text{tol}} = 10^{-6}$ in the out-of-sample period in order to guarantee strict adherence to the constraint set. Figure \ref{fig:learn_A}(a) reports the average training loss at each epoch and the $95\%$-ile confidence interval, evaluated over $10$ independent trials. Once again, we observe that the loss curves for the ADMM model and OptNet model are very similar. Interestingly, we observe that the ADMM model produces a consistently lower average training loss. Recall that both models use implicit differentiation to compute the relevant gradient, which assumes an exact fixed point at each optimal solution $\toi{\bz^*(\btheta)}$. In practice, each  $\toi{\bz^*(\btheta)}$ is only approximately optimal, to within a tolerance  $\epsilon_{\text{tol}}$, and therefore  differentiating at a solution that is not an exact fixed point will result in small errors in the gradient that likely explain the observed difference. That said, the training loss profile of the ADMM and OptNet models are very similar, and the final models achieve approximately equal loss after $500$ epochs. Figure \ref{fig:learn_A}(b) compares the average and $95\%$-ile confidence interval of the total time spent executing the forward and backward pass algorithms during training. Once again we observe that the ADMM model is shown to be approximately $5$ times faster than the OptNet model and requires less than $100$ seconds to train.

\begin{figure}[H]
  \centering
      \includegraphics[width=\linewidth,height = 5cm, trim={0cm 0cm 0cm 0cm},clip]{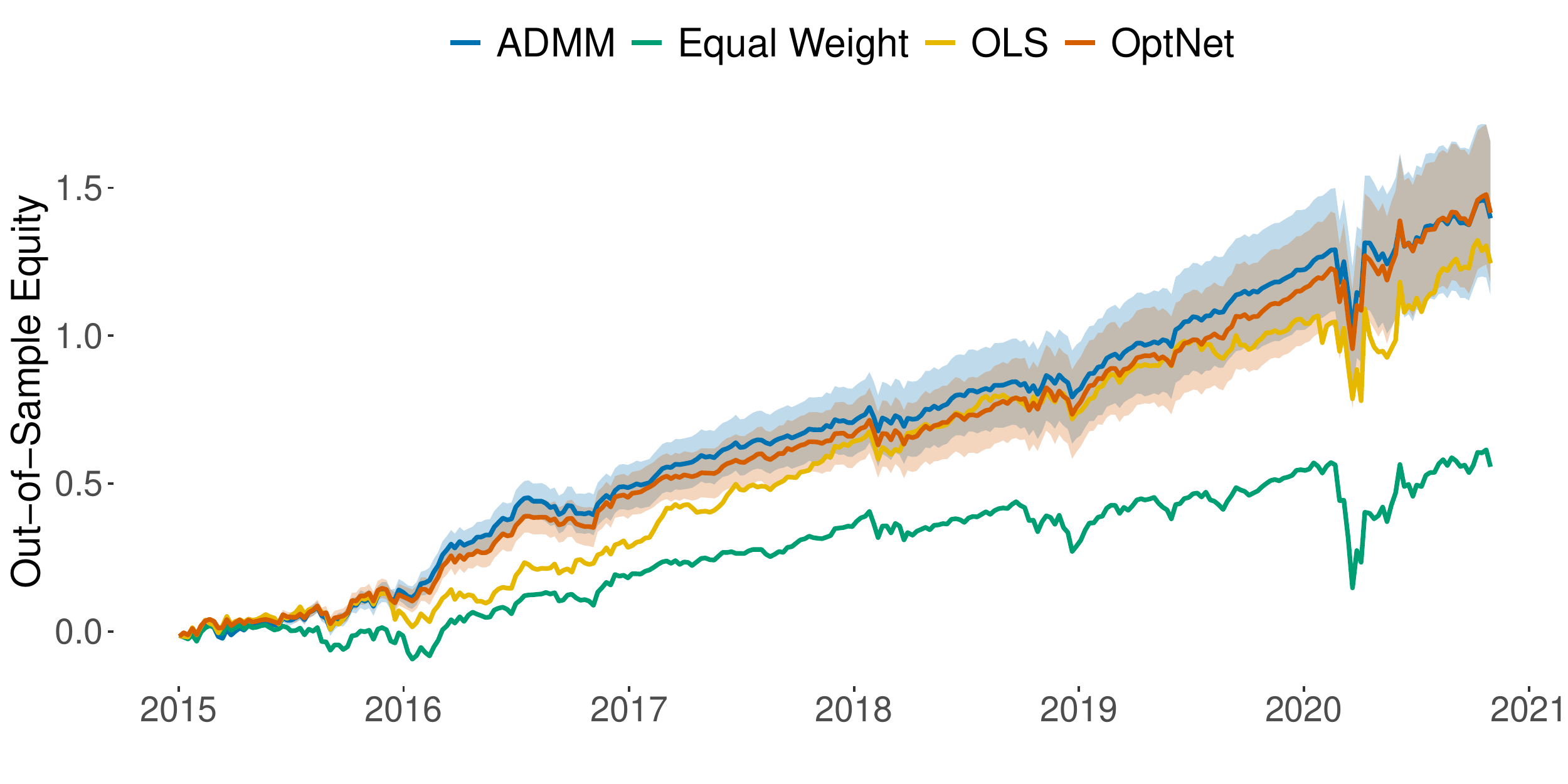}
  \caption{Out-of-sample equity growth for ADMM IPO max-Sharpe portfolio, Equal Weight portfolio, OLS max-Sharpe portfolio and OptNet IPO max-Sharpe portfolio.}
  \label{fig:eq_msr_admm}
\end{figure}

Figure \ref{fig:eq_msr_admm} reports the out-of-sample equity growth of the ADMM IPO max-Sharpe portfolio, Equal Weight portfolio, OLS max-Sharpe portfolio and OptNet IPO max-Sharpe portfolio. The out-of-sample economic performance metrics are reported in Table \ref{tab:msr_admm}.  First, observe that all max-Sharpe models outperform the Equal Weight benchmark on an absolute and risk-adjusted basis. Furthermore, the ADMM and OptNet IPO max-Sharpe models achieve an out-of-sample Sharpe ratio  that is approximately $50\%$ higher than that of the naive `predict, then optimize' OLS max-Sharpe model, thus highlighting the benefit of training a fully integrated system. Lastly, the ADMM model achieves a marginally higher out-of-sample Sharpe ratio in comparison to the OptNet model, though the difference is not statistically significant.

\begin{table}[H]
\centering
\begin{tabular}{l | l l l l}
\toprule
                 & ADMM & Equal Weight & OLS     & OptNet\\
\hline
Mean             & 0.2382 & 0.0950     & 0.2122  & 0.2413\\
Volatility       & 0.1777 & 0.1950     & 0.2435  & 0.1880  \\
Sharpe Ratio     & 1.3407 & 0.4872     & 0.8747  & 1.2836\\
\hline
\end{tabular}
\caption{Out-of-sample economic performance metrics for ADMM IPO max-Sharpe portfolio, Equal Weight portfolio, OLS max-Sharpe portfolio and OptNet IPO max-Sharpe portfolio. }
\label{tab:msr_admm}
\end{table}

\subsection{Experiment 4: learning $\bQ$ }\label{sec:results_learn_Q}
We consider another real-world experiment from portfolio optimization whereby the objective is to learn a parameterized model for the variable $\bQ$. We use the  same asset universe of $d_z = 255$ liquid US stocks from Experiment 3 and an identical experimental design. Here, we consider the long-only, fully-invested minimum variance portfolio optimization, described in Program $\eqref{eq:min_var}$.
\begin{equation}\label{eq:min_var}
\begin{split}
\minimize_{\bz} \frac{1}{2} \quad & \bz^T \toi{\bQ} \bz \\
\text{subject to }\quad  & \bone^T\bz = 1, \quad   0\leq \bz \leq 1.\\
\end{split}
\end{equation}
As before, the learning process is posed as a bi-level optimization program where  the objective is to learn a parameter $\btheta$ and the associated covariance matrix, $\toi{\hat{\bQ}(\btheta)}$, in order to minimize the average realized variance induced by the optimal decision policies $\{\toi{\bz^*(\btheta)}\}_{i=1}^m$.

\begin{equation} \label{eq:min_var_bi}
 \begin{split}
 \minimize_{ \btheta } \quad & \frac{1}{m} \sum_{i = 1}^m  \toit{\bz^*(\btheta)} \toi{\bQ} \toi{\bz^*(\btheta)} \\
 \text{subject to }\quad  &  \bz^*(\btheta)   = \argmin_{\bz}  \frac{1}{2} \bz^T \toi{\hat{\bQ}(\btheta)} \bz \quad \forall i=1,...,m\\
 \quad & \bone^T \toi{\bz^*(\btheta)} =1, \quad  0 \leq \toi{\bz^*(\btheta)}  \leq 1 \quad \forall i=1,...,m
 \end{split}
 \end{equation}

Asset returns, $\toi{ \ba }$  are modelled using the Famma-French Five (FF5) factor model. We follow \citet{Butler2021IPOb} and model $\toi{\bw}$ according to a multivariate GARCH$(1,1)$ process with  constant  correlation. We let $\toi{\hat{\bW }}$ denote the the time-varying covariance estimate of the auxiliary feature variables. We therefore model the stock covariance matrix as follows:
\begin{equation} \label{eq:y_hat}
\begin{split}
\toi{\hat{\ba}} & =   \btheta^T \toi {\bw}, \quad \toi{\hat{ \bQ } } =  \btheta^T  \toi{\hat{\bW}}  \btheta + \hat{\bF},
\end{split}
\end{equation}
where  $\hat{\bF}$ denotes the diagonal matrix of residual variances.

Again, the goal is to observe the training and out-of-sample performance of the ADMM model in comparison to the OptNet model. As a benchmark, we include the out-of-sample performance of the equal weight portfolio, and a minimum variance portfolio where $\btheta$ is fit by OLS. The training process for each trial consists of $200$ epochs with a mini-batch size of $32$. Portfolio models are fit to an accuracy of $\epsilon_{\text{tol}} = 10^{-4}$ in training, and a higher accuracy of $\epsilon_{\text{tol}} = 10^{-6}$ in the out-of-sample period.

\begin{figure}[H]
  \centering
  \begin{subfigure}[b]{0.45\linewidth}
    \includegraphics[width=\linewidth, trim={0cm 0cm 0cm 0cm},clip]{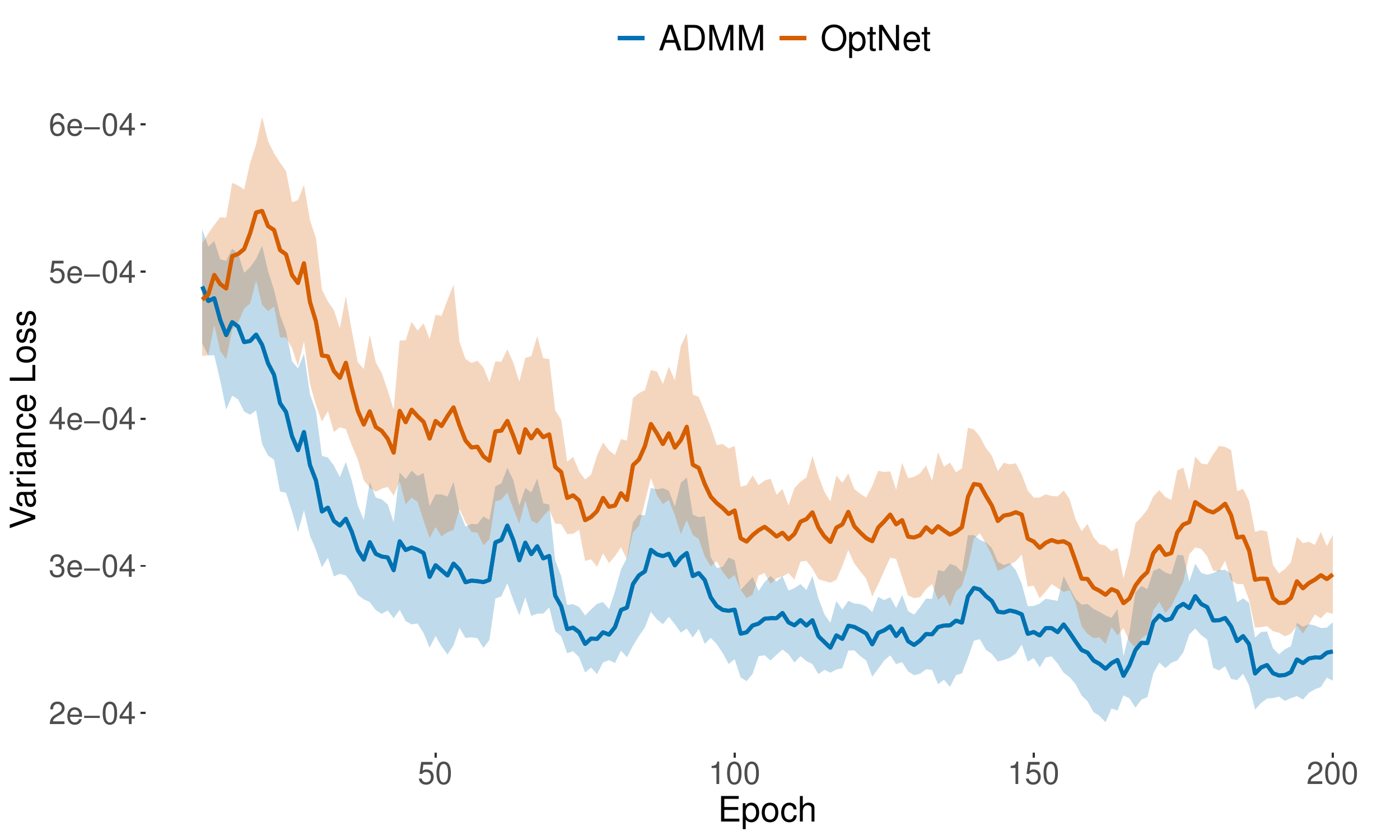}
    \caption{Training Loss}
  \end{subfigure}
  \begin{subfigure}[b]{0.45\linewidth}
    \includegraphics[width=\linewidth, trim={0cm 0cm 0cm 0cm},clip]{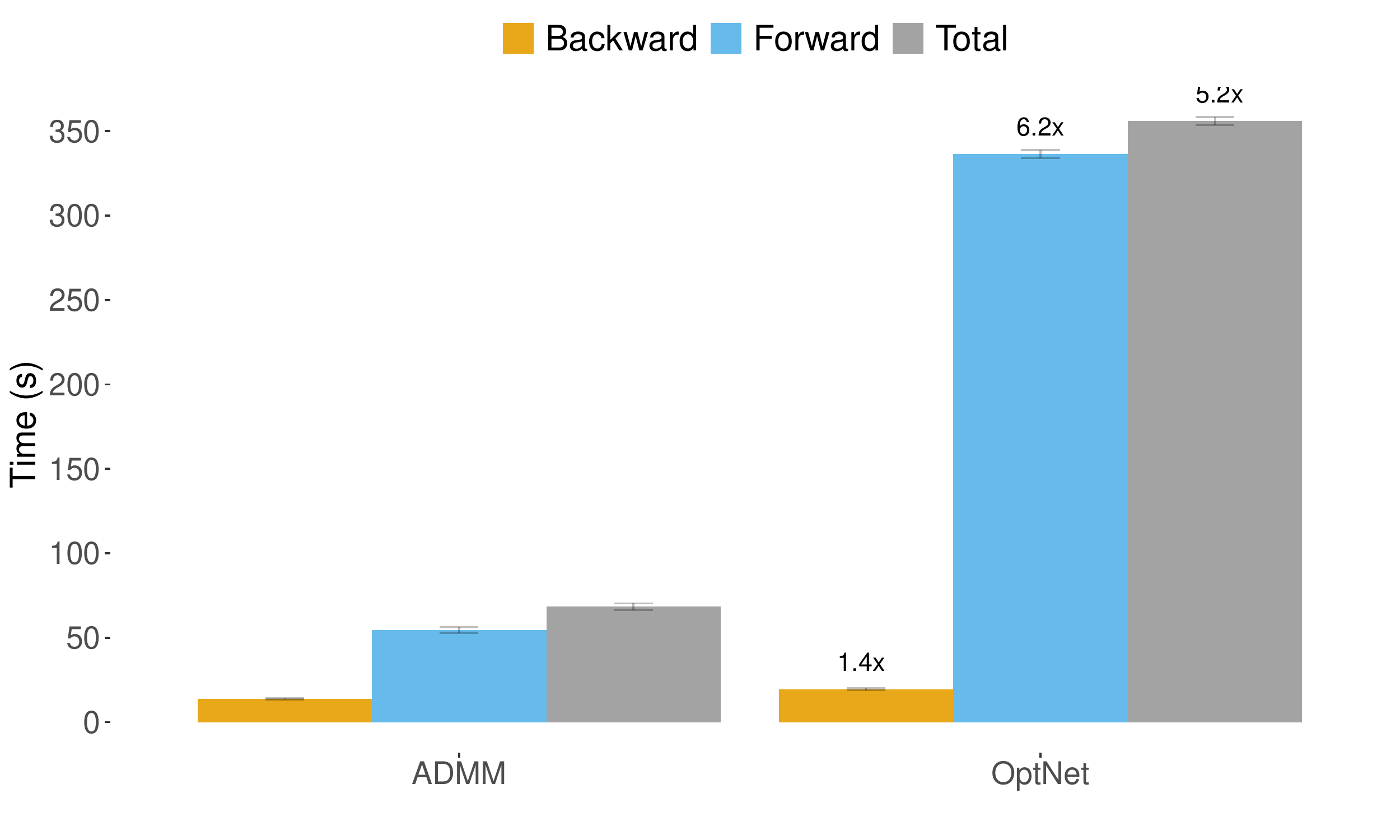}
    \caption{Computational Performance}
  \end{subfigure}
  \caption{Training loss and computational performance for learning $\bQ$ on US stock data. \\ Batch size $=32$ and $d_z = 255$.}
  \label{fig:learn_Q}
\end{figure}

Figure \ref{fig:learn_Q}(a) reports the average training loss at each epoch and the $95\%$-ile confidence interval, evaluated over $10$ independent trials. We observe that the loss curves for the ADMM model and OptNet model are very similar, thus highlighting the accuracy of the ADMM-layer. Once again we observe that the ADMM model produces a consistently lower average training loss and we refer to the discussion in Experiment 3 for a likely explanation. Figure \ref{fig:learn_Q}(b) compares the average and $95\%$-ile confidence interval of the total time spent executing the forward and backward pass algorithms in training. As before, we observe that the ADMM model requires approximately $60$ seconds to train and  is approximately $5$ times faster than the OptNet model, which requires over $300$ seconds.

\begin{figure}[H]
  \centering
      \includegraphics[width=\linewidth,height = 5cm, trim={0cm 0cm 0cm 0cm},clip]{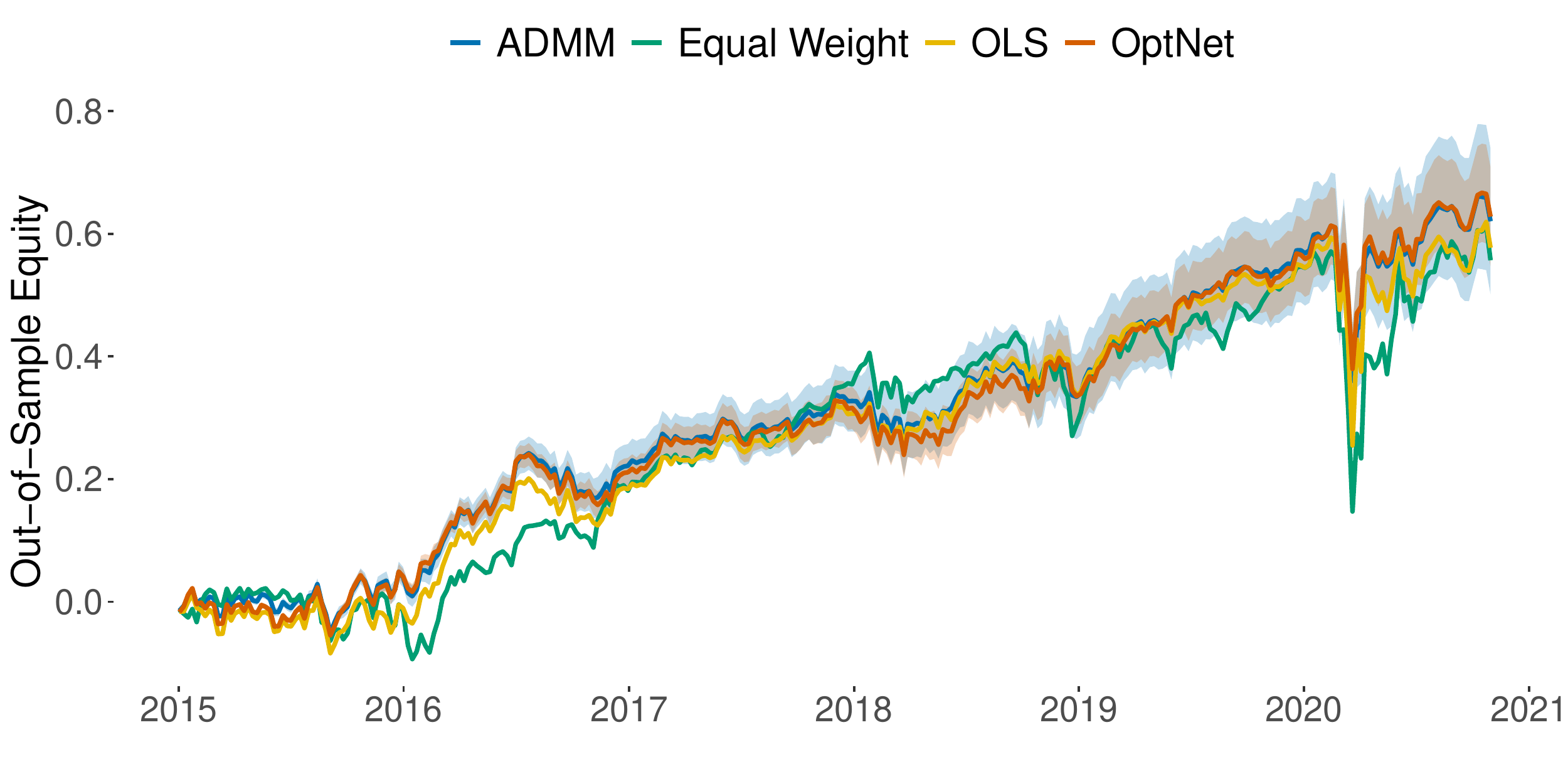}
  \caption{Out-of-sample equity growth for ADMM IPO minimum variance portfolio, Equal Weight portfolio, OLS minimum variance portfolio and OptNet IPO minimum variance portfolio.}
  \label{fig:eq_mv_admm}
\end{figure}

Finally, Figure \ref{fig:eq_mv_admm} reports the out-of-sample equity growth of the ADMM IPO minimum variance portfolio, Equal Weight portfolio, OLS minimum variance portfolio and OptNet IPO minimum variance portfolio. The out-of-sample economic performance metrics are reported in Table \ref{tab:mv_admm}.  Again, observe that all minimum-variance models outperform the Equal Weight benchmark on an absolute and risk-adjusted basis. Furthermore, the ADMM and OptNet IPO minimum variance models achieve an out-of-sample  volatility that is approximately $25\%$ lower and Sharpe ratio  that is approximately $35\%$ higher than that of the naive `predict, then optimize' OLS minimum variance model. These results are broadly consistent with the findings in  \citet{Butler2021IPOb}, who consider an identical stock universe but with considerably smaller portfolios  $(d_z \leq 100$). Our ADMM model, on the other hand, is able to overcome the computational challenges in the medium to large scale limit (described in \citet{Butler2021IPOb}), without any apparent loss in performance accuracy. 

\begin{table}[H]
\centering
\begin{tabular}{l | l l l l}
\toprule
                 & ADMM & Equal Weight & OLS     & OptNet\\
\hline
Mean             & 0.1058 & 0.0950     & 0.0984  & 0.1070\\
Volatility       & 0.1420 & 0.1950     & 0.1786  & 0.1443 \\
Sharpe Ratio     & 0.7446 & 0.4872     & 0.5510  & 0.7418\\
\hline
\end{tabular}
\caption{Out-of-sample economic performance metrics for ADMM IPO minimum variance portfolio, Equal Weight portfolio, OLS minimum variance portfolio and OptNet IPO minimum variance portfolio. }
\label{tab:mv_admm}
\end{table}

\section{Conclusion and future work}\label{sec:conclusion}
In this paper, we provide a novel and efficient framework for differentiable constrained quadratic programming. Our differentiable quadratic programming layer is built on top of the ADMM algorithm, which for medium to large scale problems is shown to be approximately an order of magnitude more efficient than the interior point implementation of the OptNet layer. The backward-pass algorithm computes the relevant problem variable gradients by implicit differentiation of a modified fixed-point iteration, which is computationally favorable to KKT implicit differentiation and memory efficient in comparison to standard unrolled differentiation. Numerical results, using both simulated and real problem data, demonstrates the efficacy of the ADMM layer, which for medium to large scale problems exhibits state-of-the art accuracy and improved computational performance.

Our experimental results should be interpreted as a proof-of-concept and we acknowledge that further testing on alternative data sets or with different problem variable assumptions is required in order to better determine the efficacy of the ADMM layer as a general purpose solver. Indeed, there is a plethora of areas for active research and algorithm improvement. First, our ADMM layer currently only supports linear equality and box inequality constraints, whereas the OptNet layer supports general linear inequality constraints. Indeed, incorporating more general inequality constraints as well as augmenting the QP with parameterized regularization norms is an active area of research.    Secondly, as discussed earlier, the ADMM algorithm is known to be vulnerable to ill-conditioned problems and the resulting convergence rates can vary significantly when the data and algorithm parameter, $\rho$, are poorly scaled. To overcome this, most first-order solvers implement a preconditioning and scaling initialization step. Currently, our ADMM-layer implementation does not support preconditioning and scaling, which is challenged by virtue of the fact that the problem data is expected to change at each epoch of training. Instead, we leave it to the user to select $\rho$ and manually scale the problem data and acknowledge that preconditioning, scaling and automatic parameter selection is an important area of future development.

Furthermore, there are several heuristic methods that could be implemented in order to improve the convergence rates and computational efficiency of our ADMM forward-pass. For example, acceleration methods, such as Andersen acceleration \citep{Anderson1965}, have recently been shown to improve the convergence rates of first-order solvers \citep{Sopa2019,Walker2011}. Indeed, applying acceleration methods to the modified fixed-point algorithm presented in this paper provides an numerically efficient scheme for potentially improving the convergence rate of the ADMM algorithm and is an interesting area of future research. Alternatively, methods that hybridize the efficiency of first-order methods with the precision of interior-point methods, such solution polishing and refinement \citep{Boyd2018}, is another area of future exploration. Nonetheless, our proposed ADMM layer is shown to be highly effective and in its current form can be instrumental for efficiently solving real-world medium and large scale learning problems.


\bibliographystyle{plainnat}
\bibliography{/Users/tars/workspace/phd_thesis/Bibliography/Bibliography}

\begin{thebibliography}{46}
\providecommand{\natexlab}[1]{#1}
\providecommand{\url}[1]{\texttt{#1}}
\expandafter\ifx\csname urlstyle\endcsname\relax
  \providecommand{\doi}[1]{doi: #1}\else
  \providecommand{\doi}{doi: \begingroup \urlstyle{rm}\Url}\fi

\bibitem[Agrawal et~al.(2019)Agrawal, Amos, Barratt, Boyd, Diamond, and
  Kolter]{Agrawal2019}
Akshay Agrawal, Brandon Amos, Shane Barratt, Stephen Boyd, Steven Diamond, and
  J.~Zico Kolter.
\newblock Differentiable convex optimization layers.
\newblock In \emph{Advances in Neural Information Processing Systems},
  volume~32, pages 9562--9574. Curran Associates, Inc., 2019.

\bibitem[Agrawal et~al.(2020)Agrawal, Barratt, Boyd, Busseti, and
  Moursi]{Agrawal2020}
Akshay Agrawal, Shane Barratt, Stephen Boyd, Enzo Busseti, and Walaa~M. Moursi.
\newblock Differentiating through a cone program, 2020.
\newblock URL \url{http://arxiv.org/abs/1703.00443}.

\bibitem[Amos and Kolter(2017)]{Amos2017}
Brandon Amos and J.~Zico Kolter.
\newblock Optnet: Differentiable optimization as a layer in neural networks.
\newblock \emph{CoRR}, abs/1703.00443, 2017.
\newblock URL \url{http://arxiv.org/abs/1703.00443}.

\bibitem[Amos et~al.(2019)Amos, Rodriguez, Sacks, Boots, and Kolter]{Amos2019}
Brandon Amos, Ivan Dario~Jimenez Rodriguez, Jacob Sacks, Byron Boots, and
  J.~Zico Kolter.
\newblock Differentiable mpc for end-to-end planning and control, 2019.

\bibitem[Anderson(1965)]{Anderson1965}
Donald G.~M. Anderson.
\newblock Iterative procedures for nonlinear integral equations.
\newblock \emph{J. ACM}, 12:\penalty0 547--560, 1965.

\bibitem[Belanger et~al.(2017)Belanger, Yang, and McCallum]{Belanger2017}
David Belanger, Bishan Yang, and Andrew McCallum.
\newblock End-to-end learning for structured prediction energy networks, 2017.

\bibitem[Bertsimas and Kallus(2020)]{Bert2020}
Dimitris Bertsimas and Nathan Kallus.
\newblock From predictive to prescriptive analytics.
\newblock \emph{Management Science}, 66\penalty0 (3):\penalty0 1025--1044,
  2020.

\bibitem[Black and Litterman(1991)]{Black1991}
F.~Black and R.~Litterman.
\newblock Asset allocation combining investor views with market equilibrium.
\newblock \emph{Journal of Fixed Income}, 1\penalty0 (2):\penalty0 7--18, 1991.

\bibitem[Blondel et~al.(2021)Blondel, Berthet, Cuturi, Frostig, Hoyer,
  Llinares-Lopez, Pedregosa, and Vert]{Blondel2021}
Mathieu Blondel, Quentin Berthet, Marco Cuturi, Roy Frostig, Stephan Hoyer,
  Felipe Llinares-Lopez, Fabian Pedregosa, and Jean-Philippe Vert.
\newblock Efficient and modular implicit differentiation, 2021.

\bibitem[Boyd and Vandenberghe(2004)]{Boyd2004}
Stephen Boyd and Lieven Vandenberghe.
\newblock \emph{Convex Optimization}.
\newblock Cambridge University Press, 2004.

\bibitem[Boyd et~al.(2011)Boyd, Parikh, Chu, Peleato, and Eckstein]{Boyd2011}
Stephen Boyd, Neal Parikh, Eric Chu, Borja Peleato, and Jonathan Eckstein.
\newblock Distributed optimization and statistical learning via the alternating
  direction method of multipliers.
\newblock \emph{Foundations and Trends in Machine Learning}, 3:\penalty0
  1--122, 01 2011.
\newblock \doi{10.1561/2200000016}.

\bibitem[Busseti et~al.(2018)Busseti, Moursi, and Boyd]{Boyd2018}
E.~Busseti, W.~Moursi, and S.~Boyd.
\newblock Solution refinement at regular points of conic problems, 2018.

\bibitem[Butler and Kwon(2021{\natexlab{a}})]{Butler2021IPOb}
Andrew Butler and Roy Kwon.
\newblock Covariance estimation for risk-based portfolio optimization: an
  integrated approach.
\newblock \emph{Social Science Research Network}, 1\penalty0 (1),
  2021{\natexlab{a}}.

\bibitem[Butler and Kwon(2021{\natexlab{b}})]{Butler2020IPO}
Andrew Butler and Roy~H. Kwon.
\newblock Integrating prediction in mean-variance portfolio optimization.
\newblock \emph{arXiv 2102.09287}, pages 1--33, 2021{\natexlab{b}}.

\bibitem[Cornuejols and Tutuncu(2006)]{Corn2006}
Gerard Cornuejols and Reha Tutuncu.
\newblock \emph{Optimization Methods in Finance}.
\newblock Cambridge University Press, 2006.

\bibitem[Diamond et~al.(2018)Diamond, Sitzmann, Heide, and
  Wetzstein]{Diamond2018}
Steven Diamond, Vincent Sitzmann, Felix Heide, and Gordon Wetzstein.
\newblock Unrolled optimization with deep priors, 2018.

\bibitem[Domke(2012)]{Domke2012}
Justin Domke.
\newblock Generic methods for optimization-based modeling.
\newblock In \emph{AISTATS}, 2012.

\bibitem[Dontchev and Rockafellar(2009)]{Dontchev2009}
Asen Dontchev and R~Rockafellar.
\newblock \emph{Implicit Functions and Solution Mappings: A View from
  Variational Analysis}.
\newblock 01 2009.
\newblock ISBN 978-0-387-87820-1.
\newblock \doi{10.1007/978-0-387-87821-8}.

\bibitem[Donti et~al.(2017)Donti, Amos, and Kolter]{Donti2017}
Priya~L. Donti, Brandon Amos, and J.~Zico Kolter.
\newblock {Task-based End-to-end Model Learning}.
\newblock \emph{CoRR}, abs/1703.04529, 2017.
\newblock URL \url{http://arxiv.org/abs/1703.04529}.

\bibitem[Elmachtoub and Grigas(2020)]{Elma2020}
Adam~N. Elmachtoub and Paul Grigas.
\newblock Smart predict, then optimize.
\newblock \emph{arXiv}, 2020.

\bibitem[Elmachtoub et~al.(2020)Elmachtoub, Liang, and McNellis]{Elma2020b}
Adam~N. Elmachtoub, Jason Cheuk~Nam Liang, and Ryan McNellis.
\newblock Decision trees for decision-making under the predict-then-optimize
  framework, 2020.

\bibitem[Fama and French(2015)]{Fama2015}
Eugene~F. Fama and Kenneth~R. French.
\newblock A five-factor asset pricing model.
\newblock \emph{Journal of Financial Economics}, 116\penalty0 (1):\penalty0 1
  -- 22, 2015.
\newblock ISSN 0304-405X.
\newblock \doi{https://doi.org/10.1016/j.jfineco.2014.10.010}.
\newblock URL
  \url{http://www.sciencedirect.com/science/article/pii/S0304405X14002323}.

\bibitem[Feng and Simon(2017)]{Feng2017}
Jean Feng and Noah Simon.
\newblock Gradient-based regularization parameter selection for problems with
  non-smooth penalty functions, 2017.

\bibitem[Gabay and Mercier(1976)]{Gabay1976}
Daniel Gabay and Bertrand Mercier.
\newblock A dual algorithm for the solution of nonlinear variational problems
  via finite element approximation.
\newblock \emph{Computers \& Mathematics With Applications}, 2:\penalty0
  17--40, 1976.

\bibitem[Ganti and Gray(2011)]{Ganti2011}
R.~Ganti and Alexander~G. Gray.
\newblock Cake: Convex adaptive kernel density estimation.
\newblock In \emph{AISTATS}, 2011.

\bibitem[Glowinski and Marroco(1975)]{Glow1975}
Roland Glowinski and A.~Marroco.
\newblock Sur l'approximation, par {\'e}l{\'e}ments finis d'ordre un, et la
  r{\'e}solution, par p{\'e}nalisation-dualit{\'e} d'une classe de
  probl{\`e}mes de dirichlet non lin{\'e}aires.
\newblock 1975.

\bibitem[Goldfarb and Liu(1991)]{Goldfarb1991}
D.~Goldfarb and Shucheng Liu.
\newblock An o(n3l) primal interior point algorithm for convex quadratic
  programming.
\newblock \emph{Mathematical Programming}, 49:\penalty0 325--340, 1991.

\bibitem[Grigas et~al.(2021)Grigas, Qi, Zuo-Jun, and Shen]{Grigas2021}
Paul Grigas, Meng Qi, Zuo-Jun, and Shen.
\newblock Integrated conditional estimation-optimization, 2021.

\bibitem[Ho et~al.(2015)Ho, Sun, and Xin]{Ho2015}
Michael Ho, Zheng Sun, and Jack Xin.
\newblock Weighted elastic net penalized mean-variance portfolio design and
  computation.
\newblock \emph{SIAM Journal on Financial Mathematics}, 6\penalty0
  (1):\penalty0 1220--1244, 2015.

\bibitem[Kim et~al.(2008)Kim, Koh, Lustig, Boyd, and Gorinevsky]{Kim2008}
Seung-Jean Kim, K.~Koh, M.~Lustig, Stephen Boyd, and Dimitry Gorinevsky.
\newblock An interior-point method for large-scale l1-regularized least
  squares.
\newblock \emph{Selected Topics in Signal Processing, IEEE Journal of},
  1:\penalty0 606 -- 617, 01 2008.
\newblock \doi{10.1109/JSTSP.2007.910971}.

\bibitem[Lorraine and Duvenaud(2018)]{Lorraine2018}
Jonathan Lorraine and David Duvenaud.
\newblock Stochastic hyperparameter optimization through hypernetworks, 2018.

\bibitem[Mandi and Guns(2020)]{Mandi2020}
Jayanta Mandi and Tias Guns.
\newblock Interior point solving for lp-based prediction+optimisation, 2020.

\bibitem[Mandi et~al.(2019)Mandi, Demirovic, Stuckey, and Guns]{Mandi2019}
Jaynta Mandi, Emir Demirovic, Peter.~J Stuckey, and Tias Guns.
\newblock Smart predict-and-optimize for hard combinatorial optimization
  problems, 2019.

\bibitem[Markowitz(1952)]{Markowitz1952}
H.~Markowitz.
\newblock Portfolio selection.
\newblock \emph{Journal of Finance}, 7\penalty0 (1):\penalty0 77--91, 1952.

\bibitem[Michaud and Michaud(2008)]{Michaud2008b}
Richard Michaud and Robert Michaud.
\newblock Estimation error and portfolio optimization: A resampling solution.
\newblock \emph{Journal of Investment Management}, 6\penalty0 (1):\penalty0
  8--28, 2008.

\bibitem[O'Donoghue et~al.(2016)O'Donoghue, Chu, Parikh, and Boyd]{Boyd2016}
Brendan O'Donoghue, Eric Chu, Neal Parikh, and Stephen Boyd.
\newblock Conic optimization via operator splitting and homogeneous self-dual
  embedding, 2016.

\bibitem[Rumelhart et~al.(1986)Rumelhart, Hinton, and Williams]{Hinton1986}
David~E. Rumelhart, Geoffrey~E. Hinton, and Ronald~J. Williams.
\newblock Learning representations by back-propagating errors.
\newblock \emph{Nature}, 323:\penalty0 533--536, 1986.

\bibitem[Schubiger et~al.(2020)Schubiger, Banjac, and Lygeros]{Schu2020}
Michel Schubiger, Goran Banjac, and John Lygeros.
\newblock Gpu acceleration of admm for large-scale quadratic programming.
\newblock \emph{Journal of Parallel and Distributed Computing}, 144:\penalty0
  55--67, 2020.
\newblock ISSN 0743-7315.
\newblock \doi{https://doi.org/10.1016/j.jpdc.2020.05.021}.
\newblock URL
  \url{https://www.sciencedirect.com/science/article/pii/S0743731520303063}.

\bibitem[Sopasakis et~al.(2019)Sopasakis, Menounou, and Patrinos]{Sopa2019}
Pantelis Sopasakis, Krina Menounou, and Panagiotis Patrinos.
\newblock Superscs: fast and accurate large-scale conic optimization.
\newblock pages 1500--1505, 06 2019.
\newblock \doi{10.23919/ECC.2019.8796286}.

\bibitem[Stellato et~al.(2020)Stellato, Banjac, Goulart, Bemporad, and
  Boyd]{Stellato2020}
Bartolomeo Stellato, Goran Banjac, Paul Goulart, Alberto Bemporad, and Stephen
  Boyd.
\newblock Osqp: an operator splitting solver for quadratic programs.
\newblock \emph{Mathematical Programming Computation}, 12\penalty0
  (4):\penalty0 637?672, Feb 2020.
\newblock ISSN 1867-2957.
\newblock \doi{10.1007/s12532-020-00179-2}.
\newblock URL \url{http://dx.doi.org/10.1007/s12532-020-00179-2}.

\bibitem[Tibshirani(1996)]{Tibshirani1996}
Robert Tibshirani.
\newblock Regression shrinkage and selection via the lasso.
\newblock \emph{Journal of the Royal Statistical Society}, 58\penalty0 (1),
  1996.

\bibitem[Tikhonov(1963)]{Tikhonov1963}
Andrey Tikhonov.
\newblock Solution of incorrectly formulated problems and the regularization
  method.
\newblock \emph{Soviet Mathematics}, \penalty0 (4):\penalty0 1035--1038, 1963.

\bibitem[Uysal et~al.(2021)Uysal, Li, and Mulvey]{Uysal2021}
Ayse~Sinem Uysal, Xiaoyue Li, and John~M. Mulvey.
\newblock End-to-end risk budgeting portfolio optimization with neural
  networks, 2021.

\bibitem[Walker and Ni(2011)]{Walker2011}
Homer Walker and Peng Ni.
\newblock Anderson acceleration for fixed-point iterations.
\newblock \emph{SIAM J. Numerical Analysis}, 49:\penalty0 1715--1735, 08 2011.
\newblock \doi{10.2307/23074353}.

\bibitem[Xie et~al.(2019)Xie, Wu, Zhong, Liu, and Lin]{Xie2019}
Xingyu Xie, Jianlong Wu, Zhisheng Zhong, Guangcan Liu, and Zhouchen Lin.
\newblock Differentiable linearized admm, 2019.

\bibitem[Yang et~al.(2017)Yang, Sun, Li, and Xu]{Yang2017}
Yan Yang, Jian Sun, Huibin Li, and Zongben Xu.
\newblock Admm-net: A deep learning approach for compressive sensing mri, 2017.

\end{thebibliography}

\appendix
\section{Proof of Proposition $\ref{prop:admm_fixed_point}$}\label{app:prop_1}
We define $\bv^k = \bx^{k+1} +\bmu^k$. We can therefore express Equation $\eqref{eq:admm_qp_z_iter_simp}$ as:
\begin{equation}\label{eq:app_z}
\bz^{k+1} =  \Pi( \bx^{k+1}  + {\bmu}^k )  = \Pi( \bv^k ),
\end{equation}
and Equation $\eqref{eq:admm_qp_mu_iter_simp}$ as:
\begin{equation}\label{eq:app_mu}
{\bmu}^{k+1} = {\bmu}^{k} +  \bx^{k+1} - \bz^{k+1} = \bv^k - \Pi( \bv^k ).
\end{equation}
Substituting Equations $\eqref{eq:app_z}$ and $\eqref{eq:app_mu}$ into Equation $\eqref{eq:admm_qp_z_iter_simp}$ gives the desired fixed-point iteration:

\begin{align}\label{eq:app_admm_fp}
\begin{bmatrix}
\bv^{k+1}\\
\betta^{k+2}
\end{bmatrix}
& =
\begin{bmatrix}
\bx^{k+2} +\bmu^{k+1}\\
\betta^{k+2}
\end{bmatrix}\\
&
= -
\begin{bmatrix}
\bQ + \rho \bI_{\bv} &   \bA^T \\
\bA & 0
\end{bmatrix}^{-1}
\begin{bmatrix}
 \bp - \rho (\bz^{k+1} - \bmu^{k+1})\\
 -\blb
\end{bmatrix}
+
\begin{bmatrix}
\bmu^{k+1}\\
0
\end{bmatrix}\\
& =
-
\begin{bmatrix}
\bQ + \rho \bI_{\bv} &   \bA^T \\
\bA & 0
\end{bmatrix}^{-1}
\begin{bmatrix}
 \bp - \rho (2\Pi(\bv^{k}) - \bv^{k})\\
 -\blb
\end{bmatrix}
+
\begin{bmatrix}
\bv^{k}\\
\betta^{k+1}
\end{bmatrix}
-
\begin{bmatrix}
\Pi(\bv^{k})\\
\betta^{k+1}
\end{bmatrix}.
\end{align}

\section{Proof of Proposition $\ref{prop:admm_grads}$}\label{app:prop_2}

We define $F \colon \mathbb{R}^{d_{v}} \times \mathbb{R}^{d_\eta} \to \mathbb{R}^{d_{v}} \times \mathbb{R}^{d_\eta}$ as:
\begin{align}\label{eq:app_F}
F(\bv,\betta)
& = -
\begin{bmatrix}
\bQ + \rho \bI_{\bv} &   \bA^T \\
\bA & 0
\end{bmatrix}^{-1}
\begin{bmatrix}
 \bp - \rho (2\Pi(\bv) - \bv)\\
 -\blb
\end{bmatrix}
+
\begin{bmatrix}
\bv\\
\betta
\end{bmatrix}
-
\begin{bmatrix}
\Pi(\bv)\\
\betta
\end{bmatrix},
\end{align}
and let
\begin{equation}\label{eq:app_M}
\bM =
\begin{bmatrix}
\bQ + \rho \bI_{\bv} &   \bA^T \\
\bA & 0
\end{bmatrix}.
\end{equation}
Therefore we have
\begin{equation}\label{eq:app_MF}
\bM F(\bv,\betta) = -
\begin{bmatrix}
 \bp - \rho (2\Pi(\bv) - \bv)\\
 -\blb
\end{bmatrix}
+
\bM
\begin{bmatrix}
\bv\\
\betta
\end{bmatrix}
-
\bM
\begin{bmatrix}
\Pi(\bv)\\
\betta
\end{bmatrix}.
\end{equation}
Taking the partial differentials of Equation $\eqref{eq:app_MF}$ with respect to the relevant problem variables therefore gives:
\begin{equation}\label{eq:app_MF}
\begin{split}
\bM \partial F(\bv,\betta)  & = -
\begin{bmatrix}
 \partial \bp \\
 -\partial \blb
\end{bmatrix}
+
\partial \bM
\begin{bmatrix}
\bv\\
\betta
\end{bmatrix}
-
\partial \bM
\begin{bmatrix}
\Pi(\bv)\\
\betta
\end{bmatrix}
-\partial \bM F(\bv,\betta)\\
&= -
\begin{bmatrix}
 \partial \bp \\
 -\partial \blb
\end{bmatrix}
-\partial \bM \Bigg[-\bM^{-1}
\begin{bmatrix}
 \bp - \rho (2\Pi(\bv) - \bv)\\
 -\blb
\end{bmatrix}
\Bigg]\\
&= -
\begin{bmatrix}
 \partial \bp \\
 -\partial \blb
\end{bmatrix}
-\partial \bM
\begin{bmatrix}
  \bx^* \\
 \betta^*
\end{bmatrix}\\
&= -
\begin{bmatrix}
  \partial \bp + \frac{1}{2} (\partial \bQ\bx^* + \partial \bQ^T\bx^*)  + \partial \bA^T\betta^*  \\
  -\partial \blb +  \partial \bA \bx^*
 \end{bmatrix}.
\end{split}
\end{equation}
 From Equation $\eqref{eq:app_MF}$ we have that the differential $\partial F(\bv,\betta)$ is given by:
 \begin{equation}\label{eq:app_partial_F}
 \partial F(\bv,\betta) = - \bM^{-1}
 \begin{bmatrix}
  \partial \bp + \frac{1}{2} (\partial \bQ\bx^* + \partial \bQ^T\bx^*)  + \partial \bA^T\betta^*  \\
  -\partial \blb +  \partial \bA \bx^*
 \end{bmatrix}.
 \end{equation}
Substituting the gradient action of Equation $\eqref{eq:app_partial_F}$ into Equation $\eqref{eq:jacob_v_etta}$ and taking the left matrix-vector product of the transposed Jacobian with the previous backward-pass gradient, $\frac{\partial \ell }{\partial \bz^*}$, gives the desired result.

\begin{equation}\label{eq:app_grads_admm}
\begin{bmatrix}
\hat{ \bd }_{\bx}  \\
\hat{ \bd }_{\betta}
\end{bmatrix}
 =
\begin{bmatrix}
\bQ + \rho \bI_{\bv} &   \bA^T \\
\bA & 0
\end{bmatrix}^{-1}
\Big[\bI_{\tilde{\bv}} - \nabla_{\tilde{\bv}} F(\tilde{\bv}(\btheta),\btheta) \Big]^{-T}
\begin{bmatrix}
D\Pi(\bv) & 0\\
0 & \bI_{\betta}
\end{bmatrix}
\begin{bmatrix}
\big( - \frac{\partial \ell }{\partial \bz^*} \big)^T \\
0
\end{bmatrix}.
\end{equation}
 From Equation $\eqref{eq:jacob_F}$ we have:
 \begin{equation}\label{eq:app_jacob_F}
 \bI_{\tilde{\bv}} - \nabla_{\tilde{\bv}} F(\tilde{\bv}(\btheta),\btheta) =
 \begin{bmatrix}
 \bQ + \rho \bI_{\bv} &   \bA^T \\
 \bA & \bzero
 \end{bmatrix}^{-1}
 \begin{bmatrix}
  - \rho (2D\Pi(\bv) - \bI_{\bv}) & 0\\
  0 & 0
  \end{bmatrix}
  +
  \begin{bmatrix}
  D\Pi(\bv) & 0\\
  0 & \bI_{\betta}
 \end{bmatrix}.
\end{equation}
Simplifying Equation $\eqref{eq:app_grads_admm}$ with Equation $\eqref{eq:app_jacob_F}$  yields the final expression:
\begin{equation}\label{eq:app_grads_admm_final}
\begin{split}
\begin{bmatrix}
\hat{ \bd }_{\bx}  \\
\hat{ \bd }_{\betta}
\end{bmatrix}
& =
 \Bigg[
\begin{bmatrix}
D\Pi(\bv) & 0\\
0 & \bI_{\betta}
\end{bmatrix}
\begin{bmatrix}
\bQ + \rho \bI_{\bv} &   \bA^T \\
\bA & 0
\end{bmatrix}
+
\begin{bmatrix}
- \rho (2D\Pi(\bv) - \bI_{\bv}) & 0\\
0 & 0
\end{bmatrix}
\Bigg]^{-1}
\begin{bmatrix}
D\Pi(\bv) & 0\\
0 & \bI_{\betta}
\end{bmatrix}
\begin{bmatrix}
\big( - \frac{\partial \ell }{\partial \bz^*} \big)^T \\
0
\end{bmatrix}.
\end{split}
\end{equation}

\section{Proof of Proposition $\ref{prop:admm_grads_box}$}\label{app:prop_3}
From the KKT system of equations $\eqref{eq:diff_sol}$ we have:
\begin{equation}\label{eq:app_d_lambda}
\bG^T\diag(\tilde{\blambda}^*) \hat{ \bd }_{\blambda} = \diag(\rho \bmu^*) \hat{ \bd }_{\blambda} = \Big( - \Big(\frac{\partial \ell }{\partial \bz^*}\Big)^T - \bQ \hat{ \bd }_{\bx} - \bA^T \hat{ \bd }_{\betta} \Big).
\end{equation}
 From Equation $\eqref{eq:kkt_partials}$ it follows that:
 \begin{equation}\label{eq:app_box_cond}
 \frac{\partial \ell   }{\partial \bl} = 0 \iff \blambda^*_- > 0 \quad \text{and} \quad \frac{\partial \ell }{\partial \bu} = 0 \iff \blambda^*_+ > 0,
\end{equation}
 and therefore Equation $\eqref{eq:app_d_lambda}$ uniquely determines the relevant non-zero gradients.
 Let $\tilde{\bmu}^*$ be as defined by Equation $\eqref{eq:mu_tilde}$, then it follows that:
 \begin{equation}\label{eq:app_d_lambda_tilde}
 \hat{ \bd }_{\blambda} = \diag(\rho \tilde{\bmu}^*)^{-1} \Big( - \Big(\frac{\partial \ell }{\partial \bz^*}\Big)^T - \bQ \hat{ \bd }_{\bx} - \bA^T \hat{ \bd }_{\betta} \Big).
 \end{equation}
Substituting $\hat{ \bd }_{\blambda}$ into Equation $\eqref{eq:kkt_partials}$ gives the desired gradients.

\section{Data Summary } \label{sec:app_imp_details}


\begin{table}[H]
\centering
\scriptsize
 \begin{tabular}{  l  c c c c c c c c }
 \hline
{\bf{GICS Sector}} & & & &  {\bf{Stock Symbols}}\\
\hline
Communication &   CBB & CMCSA & DIS & FOX & IPG & LUMN & MDP & NYT \\
Services      &    T & VOD & VZ\\
\\
Consumer      &    BBY & CBRL & CCL & F & GPC & GPS & GT & HAS \\
Discretionary &    HD & HOG & HRB & JWN & LB & LEG & LEN & LOW \\
              &    MCD & NKE & NVR & NWL & PHM & PVH & ROST & TGT\\
              &    TJX & VFC & WHR & WWW\\
\\
Consumer      &   ADM & ALCO & CAG & CASY & CHD & CL & CLX & COST\\
Staples       &   CPB & FLO & GIS & HSY & K & KMB & KO & KR \\
              &   MO & PEP & PG & SYY & TAP &  TR & TSN & UVV\\
              & WBA  & WMK & WMT\\
\\
Energy        &   AE & APA  & BKR & BP & COP & CVX & EOG & HAL\\
              &   HES & MRO & OKE & OXY & SLB & VLO & WMB & XOM\\
\\
Financials    & AFG & AFL & AIG & AJG & AON & AXP & BAC & BEN \\
              & BK & BXS &  C & GL & JPM  & L  & LNC & MMC \\
              & PGR & PNC & RJF & SCHW & STT & TROW & TRV & UNM\\
              & USB & WFC & WRB & WTM\\
\\
Health        &   ABMD & ABT & AMGN & BAX & BDX & BIO & BMY & CAH\\
Care          & CI & COO & CVS & DHR & HUM & JNJ & LLY & MDT \\
              & MRK & OMI & PFE & PKI & SYK & TFX & TMO & VTRS\\
              & WST\\
\\
Industrials   &   ABM & AIR & ALK & AME & AOS & BA & CAT & CMI\\
              &   CSL & CSX & DE  & DOV & EFX & EMR & ETN & FDX\\
              &   GD & GE & GWW & HON & IEX & ITW & JCI & KSU \\
              & LMT & LUV & MAS & MMM & NOC & NPK&   NSC & PCA\\
              & RPH & PNR & ROK & ROL & RTX & SNA & SWK & TXT\\
              &   UNP\\
\\
Information &   AAPL & ADBE & ADI & ADP & ADSK & AMAT & AMD & GLW \\
Technology  & HPQ & IBM & INTC & MSFT & MSI & MU & ORCL & ROG \\
            & SWKS & TER & TXN & TYL& WDC & XRX\\
\\
Materials  &   APD & AVY & BLL & CCK & CRS & ECL & FMC & GLT \\
          & IFF & IP &  MOS & NEM & NUE & OLN & PPG & SEE  \\
          & SHW & SON & VMC\\
\\
Real      &   ALX & FRT & GTY & HST & PEAK & PSA & VNO & WRI\\
Estate    & WY \\
\\
Utilities   & AEP &  ATO  & BKH & CMS & CNP & D & DTE & DUK\\
            & ED & EIX & ETR & EVRG & EXC & LNT & NEE & NFG \\
            & NI & NJR & OGE &  PEG & PNM & PNW & PPL & SJW\\
            & SO & SWX & UGI & WEC & XEL\\
\hline
    \end{tabular}
\caption{U.S. stock data, sorted by GICS Sector. Data provided by Quandl.}
\label{table:stocks}
\end{table}

\end{document}